\documentclass[11pt]{amsart}%
\usepackage{amsmath}%
\usepackage{amssymb}%
\usepackage{amscd}%
\usepackage{color}%
\usepackage{eucal}%
\usepackage{mathrsfs}%
\usepackage[all]{xy}%
\def\ol#1{\overline{#1}}
\def\wh#1{\widehat{#1}}
\def\wt#1{\widetilde{#1}}
\theoremstyle{plain}
    \newtheorem{theorem}{Theorem}[section]
     
    \newtheorem{proposition}[theorem]{Proposition}
    \newtheorem{lemma}[theorem]{Lemma}
    \newtheorem{corollary}[theorem]{Corollary}
\theoremstyle{definition}
    \newtheorem{definition}[theorem]{Definition}

    \newtheorem{remark}[theorem]{Remark}

\def\Alphabet{A,B,C,D,E,F,G,H,I,J,K,L,M,N,O,P,Q,R,S,T,U,V,W,X,Y,Z}
\def\alphabet{a,b,c,d,e,f,g,h,i,j,k,l,m,n,o,p,q,r,s,t,u,v,w,x,y,z}
\def\endpiece{xxx}
\def\makeAlphabet[#1]{\expandafter\makeA#1,xxx,}
\def\makealphabet[#1]{\expandafter\makea#1,xxx,}
\def\makeA#1,{\def\temp{#1}\ifx\temp\endpiece\else%
\mkbb{#1}\mkfrak{#1}\mkbf{#1}\mkcal{#1}\mkscr{#1}\expandafter\makeA\fi}%
\def\makea#1,{\def\temp{#1}\ifx\temp\endpiece\else\mkfrak{#1}\mkbf{#1}\expandafter\makea\fi}%
\def\mkbb#1{\expandafter\def\csname bb#1\endcsname{\mathbb{#1}}}
\def\mkfrak#1{\expandafter\def\csname fr#1\endcsname{\mathfrak{#1}}}
\def\mkbf#1{\expandafter\def\csname b#1\endcsname{\mathbf{#1}}}
\def\mkcal#1{\expandafter\def\csname c#1\endcsname{\mathcal{#1}}}
\def\mkscr#1{\expandafter\def\csname s#1\endcsname{\mathscr{#1}}}
\def\makeop[#1]{\xmakeop#1,xxx,}
\def\mkop#1{\expandafter\def\csname #1\endcsname{{\mathrm{#1}}}} %
\def\xmakeop#1,{\def\temp{#1}\ifx\temp\endpiece\else\mkop{#1}\expandafter\xmakeop\fi}%
\makeAlphabet[\Alphabet]
\makealphabet[\alphabet]
\makeop[Alt,End,Ext,Hom,Sym,Tor]
\makeop[CH,Pic,div]
\makeop[Ker,Im,Coker,Coim,id,pr]
\makeop[Spec,Spf,Spm,Isom]
\makeop[Re,Im]%
\makeop[dR,Nis,crys,rig,syn,tor,Cont]
\makeop[Gal,Res,ab,res]
\makeop[pol]
\makeop[Var,an,Isoc,univ,Cusp]%
\makeop[Lie,coLie,MIC,DR,Gr,Ab,Cone]
\makeop[Meas,For]
\makeop[Gl,Sl,vec]
\makeop[Eis,tr,Fil,gr,et,red,loc,Fr]

\def\pair#1{\left\langle #1 \right\rangle}

\def\Isocda{\Isoc^{\kern-0.5mm\dagger}}

\def\bsK{K}

\def\col{{\operatorname{col}}}

\def\Exp#1{\exp\left[#1\right]}%
\def\pair#1{\langle#1\rangle}%
\def\markout#1{}%

\begin{document}
\title[$p$-adic Eisenstein-Kronecker series]{$p$-adic Eisenstein-Kronecker series for CM elliptic curves and the Kronecker limit formulas}
\author[Bannai]{Kenichi Bannai}
\address{Department of Mathematics, Keio University, 3-14-1 Hiyoshi, Kouhoku-ku, Yokohama 223-8522, Japan}
\email{bannai@math.keio.ac.jp}
\author[Furusho]{Hidekazu Furusho}
\address{Graduate school of Mathematics, Nagoya University, Furo-cho Chikusa-ku, Nagoya 464-8602, Japan}
\email{furusho@math.nagoya-u.ac.jp}
\author[Kobayashi]{Shinichi Kobayashi}
\address{Mathematical Institute, Tohoku University, 6-3 Aramaki-aza-Aoba, Aoba-ku, Sendai  980-8578, Japan}
\email{shinichi@math.tohoku.ac.jp}
\thanks{This research was supported in part by KAKENHI 21674001, 26247004}
\date{\today}

\begin{abstract}
	Consider an elliptic curve defined over an imaginary quadratic field $K$ with good reduction at the primes above 
	$p\geq 5$ and has complex multiplication by the full ring of integers $\cO_K$ of $K$.
	In this paper, we construct $p$-adic analogues of the Eisenstein-Kronecker series for such elliptic curve as Coleman functions on the 
	elliptic curve.  We then prove $p$-adic analogues of the first and second Kronecker limit formulas by using the distribution 
	relation of the Kronecker theta function.
\end{abstract}

\maketitle

\setcounter{tocdepth}{1}

%
%
%
\section{Introduction}\label{section: introduction}
%
%
%

Let $\Gamma \subset \bbC$ be a lattice.   Then for an integer $a$ and $z, w \in \bbC$,
the Eisenstein-Kronecker-Lerch series for the lattice $\Gamma$ is defined by
\begin{equation} \label{equation: definition of Eisenstein-Kronecker}
	K^*_a(z,w,s; \Gamma) = \sum_{\gamma \in \Gamma \setminus\{-z\}} 
    		\frac{(\ol z+ \ol\gamma)^a}{|z+\gamma|^{2s}}\chi_{w}(\gamma),
\end{equation}
where $A(\Gamma)$ is the area of the fundamental domain of $\Gamma$ divided by $\boldsymbol\pi = 3.1415\cdots$
and $\chi_w(z) :=  \Exp {(z \overline{w}  - w \overline{z})/A(\Gamma)}$ for any $z$, $w \in \mathbb{C}$.
The above series converges for $\operatorname{Re}(s) > a/2+1$, but one may give it meaning for general $s$ by analytic continuation.
In what follows, we will omit $\Gamma$ from the notation if there is no fear of confusion.
For integers $m, n$, the classical Eisenstein-Kronecker function, more commonly known as the Kronecker double series, is defined to be the function
$$
	E_{m,n}(z) : = K^*_{n-m}(0,z,n).
$$
In the context of polylogarithms, these functions are elliptic analogues of the 
Bloch-Wigner-Zagier polylogarithm function (see \cite{L} p.282 Remark.)
In the paper \cite{Col1},  Coleman defined the $p$-adic analogue of the classical polylogarithm
function as a \textit{Coleman function}, which is a class of $p$-adic analytic functions generalizing 
the rigid analytic functions.  The purpose of this paper is to define for integers $m, n$ the
$p$-adic analogue 
$
	E^\col_{m,n}(z)
$
of the Eisenstein-Kronecker function as a Coleman function, when the complex torus $\bbC/\Gamma$
has a model as  an elliptic curve defined over an imaginary quadratic field $K$ with good reduction at the primes
above $p \geq 5$ and complex multiplication by the full ring of integers $\cO_K$ of $K$.
The main ingredient in the construction of the
Eisenstein-Kronecker series is the distribution relation.

Focusing on the distribution relation of the Eisenstein-Kronecker series, we then prove the $p$-adic 
analogues of the Kronecker limit formulas.
We let $\theta(z)$ be the reduced theta function on $\bbC/\Gamma$ associated to 
the divisor $[0] \subset \bbC/\Gamma$, normalized so that $\theta'(0) =1$ (See \eqref{eq: transformation}
for the precise transformation formula.)  
Then the Kronecker limit formulas in the classical complex case are given as follows.

\begin{theorem}[Kronecker limit formulas]\label{thm: KLF}
	Let  $c$ be the Euler constant 
	$
		c:=\lim_{n\rightarrow \infty} \left(1+\frac{1}{2}+ \cdots +\frac{1}{n}-\log n\right),
	$ 
	and let $\Delta$ be the discriminant of  $\Gamma$ defined as $\Delta:=g_2^3-27g_3^2$, where 
	$
		g_k:=\sum_{\gamma \in \Gamma \setminus \{0\}} \gamma^{-2k}.
	$
 	Then we have the following.
	\begin{enumerate}
		\item The first limit formula
		$$
			\lim_{s \rightarrow 1}\left(A K^*_0(0, 0, s)-\frac{1}{s-1}\right)
			=-\frac{1}{12}\log|\Delta|^2-2\log A+2c.
		$$
		\item For  $z \notin \Gamma$, the second limit formula
		$$
			A K^*_0(0,z,1)=-\log |\theta(z)|^2+\frac{|z|^2}{A}-\frac{1}{12}\log|\Delta|^2.
		$$	
	\end{enumerate}
\end{theorem}

Numerous proofs exist for the classical case, but many of the proofs rely on arguments concerning the 
moduli space.  We give a new proof of the above theorem, valid for a fixed lattice $\Gamma \subset \bbC$,
using the Kronecker theta function and the distribution relation.   
Our view of understanding the Kronecker limit formulas in terms of the Kronecker theta function
and the distribution relation allows us to prove the $p$-adic analogues of Theorem \ref{thm: KLF}.
Suppose now that $\Gamma$ corresponds to a period lattice corresponding to the invariant differential $\omega = dx/y$
of an elliptic curve $E: y^2 = 4x^3 - g_2 x - g_3$ with complex multiplication by
the ring of integers $\cO_\bsK$ of an imaginary quadratic field $\bsK$.  We assume in addition that $E$ is defined over $\bsK$,
and that the model above has good reduction at the primes above $p \geq 5$.
We denote
$$
	K^\col_{n-m}(0,z,n) := E^\col_{m,n}(z)
$$
to highlight the analogy.  
Then in analogy with Theorem \ref{thm: KLF} (ii), we have the following.

\begin{theorem}[$p$-adic second Kronecker limit formula]\label{thm: pKLF}
	For any prime $p \geq 5$ of good reduction, we have the second limit formula
	$$
		 K^\col_0(0,z,1)= - \log_p \theta(z) - \frac{1}{12}\log_p \Delta,
	$$
	where $\log_p \theta(z)$ is a certain $p$-adic analogue of the function $\log|\theta(z)| - |z|^2/A$ 
	defined in Definition \ref{def: log-p-theta} using the reduced theta function $\theta(z)$ and the branch of our $p$-adic logarithm. 
\end{theorem}

$p$-adic analogues of the second Kronecker limit formula were previously investigated by Katz \cite{Ka} and de Shalit \cite{dS}
in the context of $p$-adic $L$-functions when $p$ is a prime of good ordinary reduction. 
Our formulation via $p$-adic Eisenstein-Kronecker series gives a direct $p$-adic analogue, and is valid even for supersingular $p$.

Suppose now that $p \geq 5$ is a prime of good \textit{ordinary} reduction.
In this case, the prime $p$ splits as $p= \frp \frp^*$ in $\cO_K$.  We denote by $\psi$ the Hecke character of $K$ associated to $E$, 
and we let $\pi:=\psi(\frp)$.
In \cite{BK1} \S 3.1, we defined a two-variable $p$-adic measure $\mu:=\mu_{0,0}$ on $\bbZ_p \times \bbZ_p$ interpolating Eisenstein-Kronecker numbers,
or more precisely, the values $K^*_{a+b}(0,0,b)/A^a$ for $a$, $b \geq 0$.
This measure depends on the choice of a $p$-adic period $\Omega_p$ of the formal group of $E$.
We define the $p$-adic Eisenstein-Kronecker-Lerch series by
$$
	K^{(p)}_{a}(0, 0, s) := \int_{\bbZ^\times_p \times \bbZ^\times_p} \pair{x}^{s-1} \omega(y)^{a-1}\pair{y}^{a-s} d \mu(x,y)
$$
for any $s \in \bbZ_p$, where $\pair{-} : \bbZ_p^\times \rightarrow \bbC_p^\times$ is given as the composition
$\bbZ_p^\times \rightarrow 1 + p\bbZ_p \hookrightarrow \bbC_p^\times$ and
$\omega: \bbZ_p^\times \rightarrow \mu_{p-1}$ is the Teichm\"uller character,
so that $x= \omega(x) \pair{x}$ for any $x \in \bbZ_p^\times$.
Then an argument similar to the proof of Theorem \ref{thm: KLF} (i) gives the following.

\begin{proposition}[$p$-adic first Kronecker limit formula]\label{pro: pKLF}
 	Suppose $p \geq 5$ is  a prime of good ordinary reduction.  Then
	$$
		\lim_{s\rightarrow 1} K^{(p)}_{0}(0,0,s) = \Omega_p^{-1}\left( 1 - \frac{1}{p} \right) \log_p \ol\pi
	$$
	where $\Omega_p$ is the $p$-adic period corresponding to $\mu$. 
\end{proposition}

The proof of the above proposition is pararell to that of the proof of Theorem \ref{thm: KLF}(i).  However, due to the
existence of a trivial zero for the function $K^{(p)}_0(0,0,s)$ at $s=1$, the analogy with the classical case is not
perfect.  See Remark \ref{rem: not pKLF} for details.\medskip

The $p$-adic analogue of the classical polylogarithm was used to express the specializations
at $p$-power  roots of unity of the $p$-adic realization of  the polylogarithm sheaf (\cite{D}, \cite{Som}, \cite{Bes1},  \cite{Bes2}, \cite{BdJ}, \cite{Ba1}, \cite{Ba4}).
The elliptic analogue of the classical polylogarithm sheaf was first constructed by Belinson and Levin \cite{BL}.
In previous research, we studied the $p$-adic realization of the elliptic polylogarithm sheaf for CM elliptic curves (\cite{Ba3}, \cite{BKT1}, \cite{BKT2}).
As in the classical case, the $p$-adic Eisenstein-Kronecker function defined in this paper should be related to specializations 
at $p$-power torsion points of the $p$-adic realization of the elliptic polylogarithm sheaf.   We expect that this function will play a
role in future research in the formulation of the $p$-adic analogue of the elliptic Zagier conjecture formulated by Wildeshaus \cite{W2}.

\tableofcontents

%
%
%
\section{Classical Kronecker limit formulas}\label{section: one}
%
%
%

In this section, we first review the definitions of the Eisenstein-Kronecker series and the Kronecker theta function.
We then give new proofs of the first and second Kronecker limit formulas using the Kronecker theta function and the distribution relation.   
Our proof in the classical complex case will be the model for proving the $p$-adic analogue.
As in the original proof by Kronecker, we first prove the second limit formula, and then 
deduce the first limit formula from the second.

%
\subsection{Eisenstein-Kronecker series and the Kronecker theta function}
%

We fix a lattice $\Gamma$ in $\mathbb{C}$ and let $A$ be the area of the fundamental domain of $\Gamma$ divided by $\pi$.  
Let $a$ be an integer and let $z_0$, $w_0 \in \bbC$.
We denote by $K^*_a(z_0,w_0,s)$ the Eisenstein-Kronecker-Lerch series given in \eqref{equation: definition of Eisenstein-Kronecker}
of the introduction.  The following result was proved by Weil \cite{We} VIII \S 13.  See also \cite{BKT1}  Proposition 2.4
for details concerning the case $a \leq 0$.

\begin{proposition} Let $a$ be an integer.
	\begin{enumerate}
		\item The function $K_a^*(z_0,w_0,s)$ for $s$ continues meromorphically to a function on the whole $s$ plane,
		with a simple pole only at $s=1$ if $a=0$ and $w_0 \in \Gamma$.
		\item The function $K^*_a(z_0,w_0,s)$ satisfies the functional equation
		$$
			\Gamma(s) K^*_a(z_0,w_0,s) = A^{a+1-2s} \Gamma(a+1-s) K^*_a(w_0,z_0,a+1-s) \chi_{w_0}(z_0).
		$$
	\end{enumerate}
\end{proposition}

As in the introduction, we define the Eisenstein-Kronecker function, referred to more commonly as the Kronecker double series, as follows.

\begin{definition} For any integer $m$, $n$, we define the \textit{Eisenstein-Kronecker function} $E_{m,n}(z)$ to be the $\sC^\infty$-function on $\bbC \setminus \Gamma$
	defined by
	$$
		E_{m,n}(z) := K^*_{n-m}(0,z,n).
	$$
\end{definition}
This function is known to satisfy the differential equations
\begin{align*}
	\partial_z E_{m,n}(z) &= - E_{m-1,n}(z)/A, &  \partial_{\ol z} E_{m,n}(z) &= E_{m,n-1}(z)/A.
\end{align*}

\begin{remark}
	 The Eisenstein-Kronecker functions $E_{m,n}(z)$ may be used to describe the $\bbR$-Hodge realization
	 of the elliptic polylogarithm sheaf (see \cite{BKT1} Appendix for details.)
\end{remark}

%
%

We next review the definition and basic properties of the Kronecker theta function.
Denote by $\theta(z)$ the reduced theta function
associated to the divisor $[0]$ of $\bbC/\Gamma$, normalized so that $\theta'(0)=1$.  
This is the function used by Robert to define his elliptic units (\cite{Rob} \S 1.)
This function satisfies the transformation formula
\begin{equation}\label{eq: transformation}
	\theta(z + \gamma) = \alpha(\gamma) \exp\left( \frac{z \gamma }{A} + \frac{|\gamma|^2}{2A} \right) \theta(z),
\end{equation}
where $\alpha(\gamma)=-1$ if $\gamma \not\in 2\Gamma$ and $\alpha(\gamma)=1$ otherwise.
We define the Kronecker theta function as follows.

\begin{definition}[Kronecker theta function]
	We let
	$$
		\Theta(z,w) :=\frac{ \theta(z+w)}{\theta(z)\theta(w)}.
	$$
	The relation of this function to the two-variable Jacobi theta function $F_\tau(z,w)$ of Zagier \cite{Zag} is given by
	$\Theta(z,w) =  \exp(z w/A)F_\tau(z, w)$.
\end{definition}

The values $K^*_a(z,w,s)$ are defined for any $z$, $w \in \bbC$, but as a function, $K^*_a(z,w,s)$ is \textit{not} 
continuous for $z$, $w \in \Gamma$.  We let $K_a(z,w,s) := K^*_a(z,w,s)$ be the $\mathscr C^\infty$ function
defined for any $z$, $w \in \bbC \setminus\Gamma$.  We regard the function $K_a(z,w,s)$ as being undefined for $z$ or $w \in \Gamma$.
Since the function $K^*_a(z,w,s)$ for $z \in \Gamma$ is defined in \eqref{equation: definition of Eisenstein-Kronecker} by removing the 
summand with poles, we have
\begin{equation}\label{eq: limit}
		\lim_{z \rightarrow 0} \left[ K_1(z,w,1) -\frac{1}{ z} \right] = K^*_1(0,w,1).
\end{equation}

The relation between this function and the Kronecker theta function is given by
the following theorem due to Kronecker. 
\begin{theorem}[Kronecker]\label{theorem; kronecker}
	$$
		\Theta(z,w) = \exp\left[ \frac{z \overline w}{A}\right] K_1(z,w,1).
	$$
\end{theorem}	
The above theorem was originally proved in terms of Jacobi theta functions by Kronecker
using moduli arguments (See for example \cite{We}.) 
In \cite{BK1} Theorem 1.13 or \cite{BK2} Theorem 2.10, we give another proof valid for a fixed lattice $\Gamma\subset\bbC$
using the fact that both sides of the equality are reduced meromorphic theta functions associated to the 
Poincar\'e bundle on $\bbC/\Gamma \times \bbC/\Gamma$,
with the same poles and the same residue at each pole.

%
\subsection{Proof of the second limit formula.}
%

We now deduce Theorem \ref{thm: KLF} (ii) from Theorem \ref{theorem; kronecker}.
	
\begin{proposition}\label{pro: up to C}
	There exists a constant $C$ such that  
	$$
		\log |\theta(z)|^2-\frac{|z|^2}{A}=-A K^*_0(0,z,1)+C
	$$
	for any $z \notin \Gamma$.
\end{proposition} 

\begin{proof}
	By Theorem \ref{theorem; kronecker}, we have 
	$$
		 \Theta(z,w) - \frac{1}{z} =
		 \exp\left[ \frac{z \ol w}{A}\right] 
		\left( K_1(z,w,1) - \frac{1}{z} \right) + 
		\frac{1}{z}  \left(\exp\left[ \frac{z \ol w}{A}\right] -1 \right).
	$$
	Hence by \eqref{eq: limit},
	we have
	$$
		\lim_{z \rightarrow 0} \left( \Theta(z,w) - \frac{1}{z} \right)\\
				= K^*_1(0,w,1) + \frac{\overline w}{A}.
	$$
	Direct computation also shows that   
	\begin{equation*}
		\lim_{z \rightarrow 0} \left(\Theta(z,w) - \frac{1}{z} \right)= \frac{\theta'(w)}{\theta(w)}.
	\end{equation*}
	Hence we have
	$$
		K^*_1(0,w,1) + \frac{\overline w}{A} =  \frac{\theta'(w)}{\theta(w)} = \frac{\partial}{\partial w} \log \theta(w).
	$$
	In particular,  if we replace the variable $w$ in the above by $z$, then we have
	\begin{equation*}\label{equation: leff}
		\frac{\partial}{\partial z} \left( \log \theta(z)  - \frac{z \overline z}{A} \right)= K^*_1(0,z,1).
	\end{equation*}
	Therefore, if we let $\Xi(z)$ be the function 
	$$
		\Xi(z) := \log |\theta(z)|^2 - \frac{|z|^2}{A},
	$$
	then we have
	\begin{align*}
		\frac{\partial}{\partial z} \Xi(z) =  K^*_1(0,z,1), \qquad	\frac{\partial}{\partial \overline z} \Xi(z) 
		=  \overline{K^*_{1}(0,z,1)}. 
	\end{align*}
	On the other hand, one can directly show that 
	\begin{align*}
		A\frac{\partial}{\partial z} K^*_0(0,z,1) =-  {K^*_1(0,z,1)}, \qquad 
		A\frac{\partial}{\partial \overline z} K^*_0(0,z,1)=- {\overline{K^*_1(0,z,1)}}.  
	\end{align*}
	(See for example \cite{BKT1} Lemma A.1. )
	Hence  $\Xi(z)+AK^*_0(0,z,1)$ must be constant. 
\end{proof}

Our goal is to determine the constant $C$.   We use the following result, which is a type of distribution relation.
In what follows, we will write $z_n \neq 0$ for $n \geq 1$ to mean $z_n \in (\frac{1}{n}\Gamma/\Gamma) \setminus \{ 0\}$ for simplicity.

\begin{lemma}[Distribution Relation]\label{lem: C}
	We have  
	$$
		\sum_{z_n \not=0} \;K^*_0(0, z_n, 1)=-\frac{2\log n}{A},
	$$
	where the sum is over all $n$-torsion points $z_n$ of $\bbC/\Gamma$ except zero.  
\end{lemma} 

\begin{proof}
	We have 
	$$
		\sum_{z_n \in \frac{1}{n}\Gamma/\Gamma} \chi_{z_n}(\gamma)=
		\begin{cases}
			n^2 \qquad &(\gamma \in n\Gamma) \\
			0 \qquad &(\gamma \notin n\Gamma). 
		\end{cases}
	$$
	Therefore
	$$
		\frac{1}{n^2}\sum_{z_n \in \frac{1}{n}\Gamma/\Gamma} \;K^*_0(0, z_n, s)=\sum_{\gamma \in n\Gamma \setminus \{0\}} 
		\frac{1}{|\gamma|^{2s}}=\frac{1}{n^{2s}}K^*_0(0, 0, s)
	$$
	when the real part of $s$ is sufficiently large, and hence for any $s$ by analytic continuation.  
	In particular, we have
	$$
		\frac{1}{n^2}\sum_{z_n\not=0} \;K^*_0(0, z_n, s)=\left(\frac{1}{n^{2s}}-\frac{1}{n^2}\right)K^*_0(0, 0, s).
	$$
	By equation \cite{We} VIII \S 13 (31) noting that $\Gamma(1)=1$,  the residue of $K^*_0(0, 0, s)$ at $s=1$ is $1/A$.
	Hence we have 
	$$
		\frac{1}{n^2}\sum_{z_n\not=0} \;K^*_0(0, z_n, 1)=-\frac{2 \log n}{n^{2}A}
	$$
	as desired.
\end{proof} 
The above lemma shows that the constant $C$ is 
$$
	C=\frac{1}{n^2-1}\left[\sum_{z_n\not=0} \left( \log |\theta(z_n)|^2-\frac{|z_n|^2}{A}\right)-2\log n\right].
$$
We will now calculate this value explicitly in terms of $\Delta$. 

\begin{proposition}\label{pro: C}
	We have 
	$$
		\frac{1}{4} \log|\Delta'|^2=-
		\sum_{z_2\not=0}\left( \log |\theta(z_2)|^2-\frac{|z_2|^2}{A}\right) 
	$$
	where $z_2$ runs through non-trivial $2$-torsion points of $\bbC/\Gamma$ and 
	$$
		\Delta'=(e_1-e_2)^2(e_2-e_3)^2(e_3-e_1)^2
	$$ 
	for $y^2=4x^3-g_2x-g_3=4(x-e_1)(x-e_2)(x-e_3)$.
\end{proposition}

\begin{proof}
	Note that 
	$$
		(x-e_1)(x-e_2)(x-e_3)=\prod_{z_2\not=0}(x-\wp(z_2)).
	$$
	Then if $\Gamma=\bbZ \omega_1+\bbZ \omega_2$,
	we may suppose that $e_1=\wp(\omega_1/2)$, $e_2=\wp(\omega_2/2)$ and 
	$e_3=\wp((\omega_1+\omega_2)/2)$. 
	Since 
	$$
		{\theta(z+w) \theta(z-w)}{\theta(z)^{-2}\theta(w)^{-2}}=\wp(w)-\wp(z),
	$$
	we have  
	$$
		{\theta\left(\frac{\omega_1+\omega_2}{2}\right)\theta\left(\frac{\omega_1-\omega_2}{2}\right)}
		{\theta\left(\frac{\omega_1}{2}\right)^{-2}\theta\left(\frac{\omega_2}{2}\right)^{-2}}
		=e_2-e_1,
	$$
	$$
		{\theta\left(\omega_1+\frac{\omega_2}{2}\right)\theta\left(\frac{\omega_2}{2}\right)}
		{\theta\left(\frac{\omega_1+\omega_2}{2}\right)^{-2}\theta\left(\frac{\omega_1}{2}\right)^{-2}}
		=e_1-e_3,
	$$
	$$
		{\theta\left(\omega_2+\frac{\omega_1}{2}\right)\theta\left(\frac{\omega_1}{2}\right)}
		{\theta\left(\frac{\omega_1+\omega_2}{2}\right)^{-2}\theta\left(\frac{\omega_2}{2}\right)^{-2}}
		=e_2-e_3.
	$$
	Hence using the transformation formula
	$$
		\theta(z+\gamma) = \alpha(\gamma) \exp\left( \frac{z \ol \gamma}{A} + \frac{\gamma\ol\gamma}{2A}\right) \theta(z)
	$$
	of $\theta(z)$, where $\gamma$ is any element in $\Gamma$ and $\alpha(\gamma) = 1$ if $\gamma \in 2\Gamma$ and $=-1$ otherwise,
	the value $\Delta'$ is  
	$$
		\exp\left[\frac{\omega_1\overline{\omega_1}+\omega_2\overline{\omega_2}+\overline{\omega_1}\omega_2}{A}\right] 
		{\theta\left(\frac{\omega_1}{2}\right)^{-4}
		\theta\left(\frac{\omega_2}{2}\right)^{-4}\theta\left(\frac{\omega_1+\omega_2}{2}\right)^{-4}}.
	$$
	Multiplying it and its complex conjugation and taking the logarithm, 
	we obtain the formula. Note that since we take the 
	logarithm of  {\it positive real} numbers, the values do not depend on the choice of the branch of the logarithm.
\end{proof}

\begin{proof}[Proof of Theorem \ref{thm: KLF} (2)]
	Since the Ramanujan $\Delta$ is given by $\Delta=2^4 \Delta'$, we have by Lemma \ref{lem: C} and Proposition \ref{pro: C}
	$$
		C=\frac{1}{3}\left(-\frac{1}{4} \log|\Delta'|^2-2\log 2\right)=-\frac{1}{12}\log|\Delta|^2.
	$$
	Our assertion now follows from Proposition \ref{pro: up to C}.
\end{proof}

%
\subsection{Proof of  the first limit formula.} 
%

We first review the results of \cite{We} VIII \S 13.
For a fixed $z_0$, $w_0 \in \bbC$,  we let
$\theta^*(t,z_0,w_0)$ be the function
\begin{equation}\label{equation: star sum}
	\theta^*(t,z_0,w_0) = {\sum_{\gamma \in \Gamma}}^* \exp(- t|z_0 + \gamma|^2/A) \chi_{w_0}(\gamma),
\end{equation}
defined for $t > 0$, where 
$\sum^*$ denotes the sum taken over all $\gamma \in \Gamma$ other than $-z_0$ if $z_0$ is in $\Gamma$.
Furthermore, we let
$$
	I(z_0, w_0, s) := \int_{1}^\infty\theta^*(t,z_0,w_0) t^{s-1} dt,
$$
which converges for all $s \in \bbC$.
Then by \cite{We} VIII \S 13 (31) we have 
\begin{multline}\label{eq: integral expression}
	A^s \Gamma(s) K^*_0(z_0, w_0, s) = I(z_0, w_0, s)
	 - \frac{\delta_{z_0}}{s} \chi_{z_0}(w_0)\\ 
	+ I(w_0, z_0, 1-s) \chi_{z_0}(w_0)  + \frac{\delta_{w_0}}{s-1},
\end{multline}
where $\delta_{x} = 1$ if  $x \in \Gamma$ and $\delta_{x} = 0$ otherwise.   The above integral
expression gives the meromorphic continuation and the functional equation of $K^*_0(z_0, w_0, s)$.
We now prove Theorem \ref{thm: KLF} (i) using the second limit formula.

\begin{proof}[Proof of Theorem \ref{thm: KLF} (i)]
	From (\ref{eq: integral expression}), we have 
	\begin{align*}
		A^{s-1} \Gamma(s) & \left(A K^*_0(0, 0, s)-\frac{1}{s-1}\right)\\
		&= I(0, 0, s)
	 	- \frac{1}{s} + I(0, 0, 1-s)  - \frac{A^{s-1} \Gamma(s)-1}{s-1}.
	\end{align*}
	Therefore, we have
	\begin{equation}\label{eq: integral expression 1}
		\lim_{s \rightarrow 1}\left(A K^*_0(0, 0, s)-\frac{1}{s-1}\right) = I(0, 0, 1)
	 	-1 + I(0, 0, 0)
		-\log A+c,
	\end{equation}
	where $c$ is the Euler constant as before and we used the fact $\Gamma'(1)=-c$. 
	On the other hand, we have for $z_0 = 0$,  $w_0 = z \not\in \Gamma$ and $s=1$,
	\begin{equation*}
		 AK^*_0(0, z, 1) = I(0, z, 1) - 1 + I(z, 0, 0).
	\end{equation*}
	Note that $\displaystyle\lim_{z \rightarrow 0}  I(0, z, 1)= I(0, 0, 1)$, and
	for $z \neq 0$, if we let 
	$$
		I^* (z, 0, s)= I(z, 0, s)-\int_{1}^\infty \exp(-t|z|^2/A) t^{s-1} dt,
	$$
	then $\displaystyle\lim_{z \rightarrow 0}  I^*(z, 0, 0)= I(0, 0, 0)$ by the definition of the sum $\sum^*$ in \eqref{equation: star sum}.
	We have 
	\begin{align*}
		\Gamma(s)-\frac{1}{s}&=
		\int_{|z|^2/A}^\infty e^{-t} t^{s-1} dt+\int_{0}^{|z|^2/A} e^{-t} t^{s-1} dt-\frac{1}{s}\\
		&=
		\int_{|z|^2/A}^\infty e^{-t} t^{s-1} dt+\int_{0}^{|z|^2/A} (e^{-t}-1) t^{s-1} dt+
		\frac{1}{s}\left[ \left(\frac{|z|^2}{A}\right)^s-1\right].
	\end{align*}
	Taking $s \rightarrow 0$, we have 
	$$
		-c=\int_{|z|^2/A}^\infty e^{-t} t^{-1} dt+\int_{0}^{|z|^2/A} (e^{-t}-1) t^{-1} dt
		+\log \left(\frac{|z|^2}{A}\right).
	$$
	Hence 
	\begin{align*}
	 	&AK^*_0(0, z, 1) = I(0, z, 1) - 1 + I^*(z, 0, 0)+\int_{1}^\infty \exp(-t|z|^2/A) t^{-1} dt \\
		&= I(0, z, 1) - 1 + I^*(z, 0, 0)
		-c-\int_{0}^{|z|^2/A} (e^{-t}-1) t^{-1} dt-\log \left(\frac{|z|^2}{A}\right).
	\end{align*}
	Therefore 
	\begin{align*}
		 \lim_{z \rightarrow 0} 
		\left(AK^*_0(0, z, 1)+\log |z|^2 \right)= I(0, 0, 1) - 1 + I(0, 0, 0)
		-c+\log A.
	\end{align*}
	Finally, combining this with (\ref{eq: integral expression 1}) and the second limit formula, we have 
	\begin{align*}
		\lim_{s \rightarrow 1}\left(A K^*_0(0, 0, s)-\frac{1}{s-1}\right)&= 
		 \lim_{z \rightarrow 0} 
		\left(AK^*_0(0, z, 1)+\log |z|^2 \right) -2\log A+2c \\
		=&-\frac{1}{12}\log|\Delta|^2-2\log A+2c. 
	\end{align*}
	This proves our assertion.
\end{proof}

%
%
%
\section{Algebraic and $p$-adic properties of the Kronecker theta function}
%
%
%

In this section, we first recall the definition and the generating function of Eisenstein-Kronecker numbers, and then 
investigate the algebraic and $p$-adic properties of the function $F_{z_0,b}(z)$ defined in Definition \ref{definition: F}.
The $p$-adic Eisenstein-Kronecker functions will be defined in the next section as the iterated Coleman integrals of $F_{z_0,b}(z)$.

%
\subsection{Eisenstein-Kronecker numbers and its generating function}
%

We define the Eisenstein-Kronecker numbers to be the special values of Eisenstein-Kronecker-Lerch series 
(see also \cite{BK1} Definition 1.5 and \cite{BKT1}  Definition 2.3.)

\begin{definition}
	Let $z_0, w_0 \in \bbC$, and let $a$ and $b$ be integers such that $(a,b)\not=(-1,1)$ if $w_0\in\Gamma$.
	The Eisenstein-Kronecker number $e^*_{a,b}(z_0,w_0)$ is defined by $e^*_{a,b}(z_0,w_0) := K^*_{a+b}(z_0,w_0,b)$.
	As in \cite{BKT1}  Definition 2.3, we let
	$$
		e^*_{a,b}(z_0) := e^*_{a,b}(0,z_0) = K^*_{a+b}(0,z_0,b)
	$$
	for $z_0 \in \bbC$ such that $z_0\not\in\Gamma$ if $(a,b)\not=(-1,1)$.
\end{definition}

For any $z_0$, $w_0 \in \bbC$, we let 
$$
	\Theta_{z_0,w_0}(z,w) := \Exp{ -\frac{z_0 \ol w_0}{A} } \Exp{ - \frac{z \ol w_0 +w \ol z_0}{A}} \Theta(z+z_0, w+ w_0).
$$
This function is known to be the generating function of Eisenstein-Kronecker numbers as follows (see \cite{BK1} \S 1.14 Theorem 1.17).

\begin{theorem}\label{theorem: 3.2}
	We have
	$$
		\Theta_{z_0,w_0}(z,w) =  \chi_{z_0}(w_0)  \frac{\delta_{z_0}}{z}+ \frac{\delta_{w_0}}{w} 
		+ \sum_{a,b\geq 0} (-1)^{a+b} \frac{e^*_{a,b+1}(z_0,w_0)}{a!A^a} z^b w^a
	$$
	in a neighborhood of the origin, where $\delta_x = 1$ if $x \in \Gamma$ and $\delta_x=0$ otherwise.
\end{theorem}

We define the function $F_{z_0, b}(z)$ as in \cite{BKT1} Definition 4.2 as follows.

\begin{definition}\label{definition: F}
	For any $z_0 \in \bbC$ and integer $b \geq 0$, we define $F_{z_0, b}(z)$ to be the meromorphic function on $\bbC$ given by the equation
	$$
		\Theta_{z_0,0}(z,w) = \sum_{b \geq 0} F_{z_0,b}(z) w^{b-1}.
	$$
\end{definition}

The choice of $z_0$ in the definition of $F_{z_0,b}(z)$ depends only on the class of $z_0$ modulo $\Gamma$.
When $z_0 = 0$, we let $F_b(z) := F_{0,b}(z)$.  Explicit calculations show that we have $F_0(z)=1$ and 
$F_1(z) = \theta'(z)/\theta(z) = \zeta(z) - e^*_{0,2}z$. By definition $\Theta_{z_0,0}(z,w) = \exp(- w \ol z_0/A) \Theta(z+z_0,w)$.
This equality and the definition of $F_1$ gives the equality $F_{z_0,1}(z) = F_{1}(z + z_0) - \ol z_0/A$.
We will later show that in the $p$-adic case, $F_{z_0,b}$ for various $z_0$ paste together to form a Coleman function.

The formula for $\Theta_{z_0,w_0}(z,w)$ as the generating function for Eisenstein-Kronecker numbers gives the following.

\begin{proposition}[Generating function]\label{pro: F generating}
	For any $b \geq 0$, the Laurent expansion of $F_{z_0, b}(z)$ at $0$ is given by
	$$
		F_{z_0, b}(z) = \frac{\delta_{z_0,b}}{z}  + \sum_{a \geq 0} (-1)^{a+b-1} \frac{e^*_{a, b}(z_0)}{a! A^a} z^a,
	$$
	where $\delta_{x,b}=1$ if $b=1$ and $x \in \Gamma$ and is zero otherwise.
\end{proposition}

\begin{proof}
	The proposition follows from Theorem \ref{theorem: 3.2} and Definition \ref{definition: F}.
\end{proof}

Next, we assume that our complex torus has an algebraic model.  Let $F$ be a number field with a fixed embedding 
$F \hookrightarrow \bbC$, and assume  that we have an elliptic curve $E$ over $F$ defined by the Weierstrass equation
\begin{equation}\label{eq: W model}
	E: y^2 = 4 x^3 - g_2 x - g_3.
\end{equation}
We let $\Gamma$ be the period lattice of $E$ with respect to the invariant differential $\omega = dx/y$.
We have a complex uniformization $\xi : \bbC/\Gamma \cong E(\bbC)$ such that $dz$ corresponds to $\omega$.  

We next define an auxiliary function $L_n(z)$, which is useful since it is an algebraic meromorphic function defined over $F$
(Proposition \ref{pro: algebraic}).

\begin{definition}\label{def: Xi}
	Let $\Xi(z,w) := \exp(- F_1(z)w) \Theta(z,w)$.  We define the connection function $L_n(z)$ by the formula
	$$
		\Xi(z,w) := 	\sum_{n \geq 0} L_n(z) w^{n-1}.
	$$
\end{definition}

\begin{remark}
	Explicit calculation shows that the connection function for small $n$ is given by $L_0(z) = 1$, $L_1(z) = 0$
	and $L_2(z) = - \frac{1}{2}\wp(z)$.
\end{remark}

The function $L_n(z)$ is a periodic with respect to $\Gamma$, hence defines a function on $\bbC/\Gamma$ holomorphic outside the points corresponding to $\Gamma$.
The relation between $F_{z_0,b}(z)$ and $L_n(z)$ is given by the formula
\begin{equation}
	F_{z_0,b}(z) = \sum_{n=0}^b \frac{F_{z_0,1}(z)^{b-n}}{(b-n)!} L_n(z + z_0).
\end{equation}

The connection function $L_n(z)$ is algebraic in the following sense.

\begin{proposition}\label{pro: algebraic}
	The functions $L_n(z)$ correspond through the uniformization $\xi$ to rational functions $L_n$ on $E$ defined over $F$.
\end{proposition}

\begin{proof}
	See \cite{BKT1}  Proposition 1.6.
\end{proof}

Assume now that $E$ has complex multiplication by the ring of integers $\cO_\bsK$ of an imaginary quadratic field $\bsK$.
In this case, the function $\Theta_{z_0, w_0}(z,w)$ satisfies the following algebraicity result.

\begin{theorem}[Damerell]
	Suppose $z_0$, $w_0$ correspond to torsion points in $\bbC/\Gamma\cong E(\bbC)$.  Then we have
	$$
		\Theta_{z_0, w_0}(z,w) - \chi_{z_0}(w_0) z^{-1} \delta_{z_0} -  w^{-1} \delta_{w_0}\in \ol\bbQ[[z, w]],
	$$
	where $\delta_x = 1$ if $x \in \Gamma$ and $\delta_x = 0$ otherwise.
\end{theorem}

\begin{proof}
	This is a reformulation of the classical theorem of Damerell.
	See \cite{BK1} Corollary 2.4 for a proof.
\end{proof}

\begin{corollary}\label{cor: algebraicity}
	Suppose $z_0$ corresponds to a torsion point in $\bbC/\Gamma \cong E(\bbC)$.  Then we have
	$$
		F_{z_0, b}(z) - \delta_{z_0,b} z^{-1} \in  \ol\bbQ[[z]],
	$$
	where $\delta_{x,b} = 1$ if $b=1$ and $x \in \Gamma$, and $\delta_{x,b} = 0$ otherwise.
\end{corollary}


We fix an isomorphism $[\, ] : \cO_\bsK \cong \End_{\ol F}(E)$ so that $\alpha \in \cO_\bsK$ acts as $[\alpha]^* \omega = \alpha \omega$
on the invariant differential $\omega = dx/y$.
For any non-zero $\alpha \in \cO_\bsK$,  we denote by $E[\alpha]$ the subgroup of $E(\ol\bbQ)$ annihilated by $\alpha$.
The function $F_{z_0,b}(z)$ is known to satisfy the following distribution relation with respect to $E[\alpha]$. 

\begin{proposition}[Distribution Relation]\label{pro: distribution F}
	The function $F_{z_0, b}(z)$ satisfies the relation
	$$
		 \sum_{z_\alpha \in  E[\alpha]}  F_{z_0 + z_\alpha, b}(z) = \alpha \ol\alpha^{1-b} F_{\alpha z_0, b}(\alpha z) 
	$$
	for any non-zero $\alpha \in \cO_\bsK$.
\end{proposition}

\begin{proof}
	By \cite{BK1} Proposition 1.16, noting that
	$$
		\Theta_{\alpha z_0,0}(\alpha z, w/\ol\alpha; \Gamma) = \ol\alpha\Theta_{N(\alpha) z_0,0}(N(\alpha)z, w;  \ol\alpha \Gamma),
	$$	
	we see that the Kronecker theta function $\Theta_{z_0,0}(z,w)$ satisfies the distribution relation
	$$
		\sum_{z_\alpha \in E[\alpha]}\Theta_{z_0+z_\alpha,0}(z,w) = \alpha \Theta_{\alpha z_0,0}(\alpha z, w/\ol\alpha)
	$$	
	for any non-zero $\alpha \in \cO_\bsK$.  Our assertion follows from the definition of $F_{z_0,b}(z)$.
\end{proof}

%
\subsection{The $p$-adic properties of $F_{z_0,b}(z)$}\label{subsection: distribution}
%

We next review the $p$-adic properties of the function $F_{z_0,b}(z)$.  We let $E$ be an elliptic curve with 
complex multiplication by the ring of integers $\cO_\bsK$ of an imaginary quadratic field $\bsK$.
We assume in addition that $E$ is defined over $\bsK$ and has good reduction at the primes above $p \geq 5$.
This implies that $p$ does not ramify in $\cO_\bsK$.
We fix a Weierstrass model of $E$ over $\cO_\bsK$ with good reduction above $p$,  and by abuse of notations,
we denote again by $E$ this model defined over $\cO_\bsK$.
Let $t = -2x/y$ be the formal parameter of $E$ at the origin, and we denote by $\wh E$ the formal group of 
$E$ with respect to the parameter $t$.  We denote by $\lambda(t)$ the formal logarithm.

Let $z_0$ be a torsion point in $E(\ol\bbQ)$, and we denote by $\wh F_{z_0,b}(t) := F_{z_0,b}(z)|_{z =\lambda(t)}$ 
the formal composition of the Laurent expansion of 
$F_{z_0,b}(z)$  at the origin with the formal power series $z = \lambda(t)$.  Note that by definition, we have
\begin{equation}\label{eq: F and L}
	\wh F_{z_0,b}(t) = \sum_{n=0}^b \frac{\wh F_{z_0,1}(t)^{b-n}}{(b-n)!} \wh L_{z_0, n}(t),
\end{equation}
where $\wh L_{z_0,n}(t)  := L_n(z + z_0)|_{z = \lambda(t)}$.   

Let $\psi:=\psi_{E/\bsK}$ be the Hecke character of $\bsK$ associated to $E$.
We let $\frp$ be a prime in $\cO_\bsK$ lifting $p$, and let $\pi := \psi(\frp)$.  
Note that if $p$ is ordinary, then $\pi$ is an element such that $p = \pi \pi^*$ with $\pi^*$ a unit in $\bsK_\frp$, and 
$\pi = -p$ if $p$ is supersingular.    
We fix an embedding $\ol\bsK \hookrightarrow \bbC_p$ such that $\pi$ maps to an element of $p$-adic absolute value $<1$ in $\bbC_p$.
By Corollary \ref{cor: algebraicity}, we may view $\wh F_{z_0,b}(t)$ as a power series with coefficients in $\bbC_p$
through this embedding.  We next review the $p$-adic properties of the power series $\wh F_{z_0,b}(t)$ through this embedding.

\begin{proposition}\label{pro: radius convergence}
	Let $z_0$ be a torsion point in $E(\ol\bbQ)$ of order prime to $\frp$.
	Then the radius of convergence of  the holomorphic part 
	$$
		\wh F_{z_0,b}(t) - \delta_{z_0,b} t^{-1}  \in \bbC_p[[t]]
	$$ of $\wh F_{z_0,b}(t)$ is one.  In other words, this
	power series defines a function on $\bbB^- := \{ t \in \bbC_p \mid |t|_p < 1 \}$ if $b \not=1$ or $z_0 \not=0$,
	and $\wh F_{1}(t) := \wh F_{0,1}(t)$ defines a function on $\bbA(0) := \{ t \in \bbC_p \mid  0 < |t|_p < 1 \}$.
\end{proposition}

\begin{proof}
	See \cite{BKT1} Proposition 4.7 for the proof.
\end{proof}

In the next section, we will prove that the power series $\wh F_{z_0,b}(t)$ paste together to
form a Coleman function on the elliptic curve minus the identity.  We next review a formula for translation by $\pi^n$-torsion points.
Let $\wh E[\pi^n] \subset \wh E(\frm_p)$ denote the group of $\pi^n$ torsion points of the formal group $\wh E$,
where $\frm_p$ is the prime ideal of $\cO_{\bbC_p}$.

\begin{lemma}[Translation]\label{lem: translation}
	Suppose $z_0$ is a torsion point in $E(\ol\bbQ)$ of order prime to $\frp$.  
	Let $t_n \in \wh E[\pi^n]$, and let $z_n$ be the image of $t_n$ through the inclusion 
	$\wh E[\pi^n] \subset E(\ol\bbQ) \subset \bbC/\Gamma$
	Then we have
	$$
		\wh F_{z_0, b}(t \oplus t_n) = \wh F_{z_0 + z_n, b}(t),
	$$
	where $\oplus$ is the formal addition law of the formal group $\wh E$.
\end{lemma}

\begin{proof}
	See \cite{BKT1} Lemma 4.13 for the proof.
\end{proof}

The above lemma gives the following corollary.

\begin{corollary}[Generating function]\label{cor: generating function}
	Suppose $z_0$ is a non-zero torsion point in $E(\ol\bbQ)$ of order prime to $\frp$.  Then for any integer $a \geq 0$, we have
	$$
		(\partial^a_{t, \log} \wh F_{z_0,b}(t))|_{t = t_n} = (-1)^{a+b-1} e^*_{a,b}(z_0 + z_n)/A^a,
	$$
	where $\partial_{t,\log}$ is the differential operator $\lambda'(t)^{-1} \partial_t $.
\end{corollary}

\begin{proof}
	Since  $\partial_{t,\log}$ is invariant under translation of the formal group,  we have by Lemma \ref{lem: translation}
	$$
	(\partial^a_{t, \log} \wh F_{z_0,b}(t))|_{t = t_n} = (\partial^a_{t, \log} \wh F_{z_0,b}(t \oplus t_n))|_{t = 0}
	=  (\partial^a_{t, \log} \wh F_{z_0 + z_n,b}(t))|_{t = 0}.
	$$
	Note that if we let $z = \lambda(t)$, then we have $\partial_z=\partial_{t,\log}$.
	Hence we have
	$$ 
		(\partial^a_{t, \log} \wh F_{z_0 + z_n,b}(t))|_{t = 0} = (\partial^a_z F_{z_0 + z_n, b}(z))|_{z=0}.
	$$
	Our assertion follows from the generating function property of $F_{z_0,b}$
	given in Proposition \ref{pro: F generating}.
\end{proof}
%
%

We let $F^{(p)}_{1}(z)$ be the function
\begin{equation}\label{eq: definition Fp}
	 F^{(p)}_{1}(z)  := F_{1}(z) - \ol\pi^{-1} F_{1}(\pi z), 
\end{equation}
which is an elliptic function corresponding to a rational function defined over $\bsK$.
Then we have $F^{(p)}_{1}(z + z_0)=F_{z_0, 1}(z) - \ol\pi^{-1} F_{\pi z_0,1}(\pi z)$,
hence
\begin{equation}\label{eq: property Fp}
	F^{(p)}_{1}(z + z_0)|_{z = \lambda(t)} = \wh F_{z_0, 1}(t) - \ol\pi^{-1} \wh F_{\pi z_0, 1}([\pi] t).
\end{equation}

%
%
\section{Construction of the $p$-adic Eisenstein-Kronecker functions}\label{section: p-adic}
%
%

Coleman integration is a theory of $p$-adic integration first developed by Coleman \cite{Col1} to
define the $p$-adic polylogarithm function.
In this section, we first review the theory of Coleman integration for curves, following the description of 
Besser (\cite{Bes2}, \cite{Bes3}) using notations coming from the theory of rigid cohomology.
We then prove that the function $F_{z_0,b}(z)$ of the previous section defines a Coleman function $F^\col_b(z)$
on the CM elliptic curve.  We then define the $p$-adic Eisenstein-Kronecker functions $E^\col_{m,n}(z)$ 
to be the iterated Coleman integral of $F^\col_b(z)$ satisfying the distribution relation.

%
\subsection{Review of Coleman integration}\label{subsection: Coleman}
%

Let $L$ be a complete subfield of $\bbC_p$, with ring of integers $\cO_L$ and residue field $k$.
Let $X$ be a smooth projective irreducible curve over $\cO_L$.   Let $U \subset X$ be an affine open 
subscheme of $X$ such that the complement $X \setminus U$ is a divisor flat over $\cO_L$.

Denote by $X^\an_{\bbC_p}$ be the rigid analytic space associated to the scheme $X_{\bbC_p}$.
Its points consist of the points of $X(\bbC_p)$.
We have the specialization morphism 
$$
	\operatorname{sp} : X_{\bbC_p}^\an \rightarrow X_{\ol\bbF_p},
$$ 
and we denote the inverse image of a point $x \in X(\ol\bbF_p)$ by $]x[$, which we call the residue disc.
The set $]x[ \subset X^\an_{\bbC_p}$ is an admissible open set of $X^\an_{\bbC_p}$. We have a set-theoretic
decomposition
$$
	X^\an_{\bbC_p} = \coprod_{x \in X(\ol\bbF_p)} ]x[,
$$
which is \textit{not} an admissible covering for the rigid topology.  

We let $\cU$ be the formal completion of $U$ with respect to the special fiber and we let $\cU_{L}$
be the rigid analytic space over $L$ associated to the formal scheme $\cU$, and we let  $\cU_{\bbC_p} := \cU_L \otimes_L \bbC_p$.
Denote by $j: \cU_{\bbC_p} \hookrightarrow X^\an_{\bbC_p}$ the natural inclusion, and we let $j^{\dagger} \cO_{\cU_{\bbC_p}}$
be the ring of functions on $\cU_{\bbC_p}$ overconvergent along $X^\an_{\bbC_p} \setminus \cU_{\bbC_p}$  (\cite{Ber} 2.1.1.3).
\begin{definition}
	For any $x \in X(\ol\bbF_p)$, we let $A(]x[)$ be the ring of functions defined by
	$$%
		 A(]x[)  := \Gamma\bigl(]x[, j^{\dagger} \cO_{\cU_{\bbC_p}}\bigr). 
	$$
\end{definition}

Since $X$ is smooth, each residue disc is isomorphic to the open disc 
$$\bbB^{-} := \{ t \in \bbC_p \mid |t|_p<1 \}$$ through a
choice of local parameter $t_x$ of $X$ at $x$.  Then we have
\begin{align*}
	 A(]x[) &  \cong \Gamma(\bbB^-, \cO_{\bbB^-}) &   &x\in U(\ol\bbF_p) \\
	A(]x[) &  \cong \bigcup_{0<r<1} \Gamma(\bbA(r), \cO_{\bbA(r)}) &   & \text{otherwise},
\end{align*}
where $\bbA(r)$ is defined to be the admissible open set of $\bbB^{-}$ defined as the 
annulus $\bbA(r) := \{ t \in \bbC_p \mid r<|t|_p<1 \}$ for any real number $r$ such that $0 < r < 1$.
Note that $A(]x[)$ is isomorphic to the ring consisting of formal power series $f(t_x) = \sum_{n \geq 0} a_n t^n_x$
which converge on $\bbB^-$ if $x \in  U(\ol\bbF_p)$, and formal power series $f(t_x) = \sum_{n = -\infty}^\infty a_n t^n_x$
for $a_n \in \bbC_p$
which converge on $\bbA(r)$ for some $r < 1$ if $x \in (X\setminus U)(\ol\bbF_p)$.
This description is independent of the choice of the parameter $t_x$.

\begin{definition}
	A branch of the $p$-adic logarithm is any locally analytic group homomorphism $\log_p : \bbC_p^\times \rightarrow
	\bbC_p$, defined by the power series
	$$
		\log_p(x) = - \sum_{n > 0} \frac{(1-x)^n}{n}
	$$ 
	for $x$ in a neighborhood of $1$.  It is characterized by the value $\log_p(p)$.
\end{definition}

Suppose a branch of the $p$-adic logarithm has been chosen.  One defines $A_{\log}(]x[)$ to be $A(]x[)$ if
$x \in U(\ol\bbF_p)$ and to be the polynomial ring in the function $\log_p(t_x)$ over $A(]x[)$ if 
$x \in (X \setminus U)(\ol\bbF_p)$.  This definition is independent up to isomorphism of the choice of the local parameter $t_x$.
Set $\Omega^1_{\log}(]x[) := A_{\log}(]x[) d t_x$.  Then one defines the ring of locally analytic functions and
one forms on $U$ by
\begin{align*}
	A_\loc(U) &:= \prod_{x \in X(\ol\bbF_p)} A_{\log}(]x[),   &
	\Omega^1_\loc(U) &:= \prod_{x \in X(\ol\bbF_p)} \Omega^1_{\log}(]x[).  
\end{align*}
We have a differential $d : A_\loc(U) \rightarrow \Omega^1_\loc(U)$, which is surjective.

Suppose $k = \bbF_q$ for $q = p^h$, and let $\cX$ and $\cU$ be the formal completions
of $X$ and $U$ with respect to the special fiber.  For simplicity, we assume that there exists a Frobenius morphism 
$\phi: \cX \rightarrow \cX$, which is a $\cO_K$-linear morphism lifting the $h$-th power $\Fr^h$ of the absolute Frobenius 
$\Fr$ on $X_k := X \otimes k$ and such that $\phi(\cU) \subset \cU$.
Then this map induces a $\bbC_p$-morphism $\phi: X^\an_{\bbC_p} \rightarrow X^\an_{\bbC_p}$
by extension of scalars.

Coleman constructs a certain subring $M(U)$ of $A_\loc(U)$, 
which we call the ring of \textit{Coleman functions} on $U$,
equipped with an integration map.   
$M(U)$ is defined so that it contains rational functions on $X$ which are regular on $U$, as well as
overconvergent functions on $\cU_{\bbC_p} \subset X^\an_{\bbC_p}$, where $\cU_{\bbC_p}$ 
is the rigid analytic space associated to $\cU$.
If we denote by $M(U)/\bbC_p$ the ring $M(U)$ modulo addition by constants,
then the integration map is a vector space map $\int : M(U) \otimes_{A(U)} \Omega^1(U) \rightarrow
M(U)/\bbC_p$ characterized by the following three properties.
\begin{enumerate}
	\item We have $d (\int \omega) = \omega$ (primitive function).
	\item We have $\int(\phi^*\omega) \equiv \phi^* (\int\omega)$  in $M(U)/\bbC_p$ (Frobenius invariance).
	\item If $g \in M(U)$, then $\int dg \equiv g$ in $M(U)/\bbC_p$.
\end{enumerate}

The construction of Coleman functions gives the following lemma.

\begin{lemma}\label{lem: criteria coleman}
	Suppose $f$ is a function in $A_\loc(U)$, and suppose $P(x)$ is a polynomial in $\bbC_p$ 
	whose roots do not contain the roots of unity.
	If we have $df \in M(U) \otimes \Omega^1(U)$ and $P(\phi^*) f \in M(U)$, then
	we have $f \in M(U)$.
\end{lemma}

By abuse of notation, for any $\omega \in  M(U) \otimes_{A(U)} \Omega^1(U)$, 
we denote by $\int \omega$ any Coleman function $F_\omega \in M(U)$ satisfying $d F_\omega = \omega$.
By construction, the Coleman function $\int \omega \in M(U)$ is determined up to addition by a constant.

It is known that the above theory of integration is independent of the choice of the branch of the $p$-adic logarithm (\cite{F} Proposition 2.3.)
Other important properties of Coleman functions are the uniqueness principle (\cite{Col1} Chapter IV and \cite{F} Proposition 2.4) and the
functorial property with respect to the morphisms of the pair $(U,X)$ (see \cite{F} Proposition 2.5.)

%
\subsection{$F_{z_0,b}(z)$ as a Coleman function}\label{subsection: F}
%

We will next show that the functions $F_{z_0, b}(z)$ defined in the previous section defines a Coleman 
function on our elliptic curve.   It is striking that the functions modified for each $z_0$ nicely paste together
to form a single Coleman function on the elliptic curve.  We fix once and for all a branch of the $p$-adic
logarithm.

Let the notations be as in \S \ref{subsection: distribution}.  In particular, 
we let $E$ be the model over $\cO_\bsK$  of our CM elliptic curve, with good reduction at the primes above $p$.
Let $L$ be a finite extension of $\bsK_\frp$ in $\bbC_p$, and by abuse of notations, we denote again by $E$ the extension
of $E$ to the ring of integers $\cO_L$ of $L$.  Then for $\pi := \psi(\frp)$, multiplication by $[\pi]$ induces a Frobenius
$
	\phi: E \rightarrow E.
$
We denote by $E(\bbC_p) :=E^\an_{\bbC_p}$ the extension to $\bbC_p$ of the rigid analytic space $E_L^\an$.
From now until the end of this paper, $z$ will denote a variable on $E(\bbC_p)$.


The reside discs of $E(\bbC_p)$ is parameterized by torsion points $z_0$ of $E$ of order prime to $\frp$.
If we let $t = -2x/y$, then this gives a local parameter of 
$E$ at the origin. Then the open disc $\{ t \in \bbC_p \mid |t|_p<1 \}$ represents the residue disc 
$]0[ \subset E(\bbC_p)$ containing the identity of $E$.
If we let $z_0$ be a torsion point of $E$ of order prime to $\frp$, and if we let
$
	]z_0[ := \tau_{z_0}(]0[)
$
for the translation $\tau_{z_0} : E \rightarrow E$ defined by $\tau_{z_0}(z) := z + z_0$, then $]z_0[$ is precisely
the residue disc containing $z_0$.  The parameter $t$ of $]0[$ then gives via translation a parameter of $]z_0[$
which gives a homeomorphism $]z_0[ \cong \bbB^-$ between the residue disc containing $z_0$ and the unit ball.
In what follows, we denote by again by $t$ the parameter of $]z_0[$ obtained in this fashion.

We first start by investigating the function $F_{z_0,1}(z)$.  We let $U := E \setminus \{0\}$.
Note that by Proposition \ref{pro: radius convergence},
the power series $\wh F_{z_0,1}(t)$ defines a function in $A(]z_0[)$ through the identification
$]z_0[ \cong \bbB^-$.  

\begin{lemma}\label{lem: interpolation}
	We let $F_1^\col$ be the function in $A_{\loc}(U)$ defined by
	$$
		F_1^\col|_{]z_0[} := \wh F_{z_0,1}(t) \,\, \in \,\,  A_{\log}(]z_0[)
	$$
	on each residue disc $]z_0[$, where $z_0$ is a torsion point in $E(\ol\bbQ)$ of order prime to $\frp$,
	including the case $z_0=0$.  
	Then $F_1^\col$ is a Coleman function on $U$.
\end{lemma}

\begin{proof}
	The differential form $d F_1 = \eta + e^*_{0,2} \omega$ is known to be a meromorphic differential form of the second kind on 
	$U$ defined over $\cO_\bsK$.  From the definition of $F_{z_0,1}(z)$ given in Definition \ref{definition: F}, we have 
	$F_{z_0,1}(z) = F_1(z+z_0) - \ol z_0/A$, hence $d F_{z_0,1} = \tau_{z_0}^* d F_1$.  Hence the calculation of the differential 
	on each residue disc gives the equality $d F_1^\col = d F_1$.
	For any Coleman function $f(z)$ on $U$, denote by $f(\pi z)$ the Coleman function $(\phi^* f)(z) =[\pi]^* f(z)$.
	Consider the function
	\begin{equation}\label{eq: F tilde}
		 \wt F^{(p)}_1(z) := \left( 1 - \frac{\phi^*}{\ol\pi} \right)F_1^\col(z)  = F_1^\col(z) - \frac{1}{\ol\pi} F_1^\col(\pi z).
	\end{equation}
	For any $z \in ]z_0[$, we have $\wt z:= \pi z \in ]\pi z_0[$.  Hence
	from the definition and the fact that $[\pi] \circ \tau_{z_0} = \tau_{\pi z_0} \circ [\pi]$, 
	we have  
	$$
		F_1^\col (\pi z)|_{]z_0[} = F_1^\col (\wt z)|_{]\pi z_0[} = \wh F_{\pi z_0,1}(\wt t)
		=\wh F_{\pi z_0,1}([\pi] t),
	$$
	where $t$ and $\wt t$ are element in $\bbB^-$ corresponding to $z$ and $\wt z$.
	Substituting the above equality into the definition \eqref{eq: F tilde} of $\wt F_1^{(p)}(z)$,
	we obtain by \eqref{eq: property Fp} the equality
	$$
		\wt F^{(p)}_1|_{]z_0[}  =  F^{(p)}_{1}|_{]z_0[},
	$$
	where $F^{(p)}_1$ is the rational function on $E$ defined in \eqref{eq: definition Fp}.  
	Hence $F^\col_1$ is an element in $A_{\loc}(U)$ such that 
	$d F_1^\col = d F_1$ and $ \left( 1 - \phi^*/\ol\pi \right) F^\col_1=F^{(p)}_1$.
	Since $F^{(p)}_1$ is a rational function and $d F_1$ is a meromorphic differential form which are both regular on $U$,
	we have in particular $F^{(p)}_1 \in M(U)$ and $d F_1 \in M(U) \otimes \Omega^1(U)$. 
	Hence by Lemma \ref{lem: criteria coleman}, the function $F^\col_1$ is in fact a Coleman function on $U$.
\end{proof}

Note that if $z_0 \neq 0$, then we have $F_1^\col|_{]z_0[}  \in A(]z_0[) \subset A_{\log}(]z_0[)$.
This shows that the value $F^\col_1(z)$ for $z \in E(\bbC_p) \setminus ]0[$ is independent from the
choice of the branch of the $p$-adic logarithm.

Since $L_n$ is a rational function on $E$ with poles only at $\{0\}$ in $E$, it is in particular a 
Coleman function on $U$.  The set of Coleman functions is a ring, and we define $F^\col_b$ as follows.

\begin{definition}
	We let $F^\col_1$ to be the Coleman function of Lemma \ref{lem: interpolation}.
	For any integer $b \geq 0$,  we define $F^\col_b$ to be the Coleman function
	\begin{equation}\label{eq: F coleman}
		F^\col_b := \sum_{n=0}^b \frac{(F^\col_1)^{b-n}}{(b-n)!} L_n
	\end{equation}
	on $U$.
\end{definition}

The functions $F^\col_b$ interpolate the power series $F_{z_0,b}(z)$ of the previous section.

\begin{proposition}[Interpolation]\label{pro: interpolation}
	For any integer $b \geq 0$, the function $F^\col_b$ on the residue disc $]z_0[$ is given by
	$$
		F^\col_b|_{]z_0[} = \wh F_{z_0,b}(t)  \in A_{\log}(]z_0[).
	$$	
\end{proposition}

\begin{proof}
	The case for $b=1$ follows from the definition of $F^\col_1$.  The case for $b>1$ follows from this case,
	noting that $L_n|_{]z_0[} = \wh L_{z_0,n}(t)$ and
	comparing the definitions of \eqref{eq: F and L} and \eqref{eq: F coleman}.
\end{proof}

\begin{proposition}[Distribution Relation]\label{pro: distribution F coleman}
	For any non-zero $\alpha \in \cO_\bsK$, the Coleman function $F^\col_b$ for any integer $b \geq 0$ satisfies the distribution relation
	\begin{equation}\label{equation: Coleman distribution}
		\sum_{z_\alpha \in  E[\alpha]}  F^\col_b(z+z_\alpha) =  \alpha \ol\alpha^{1-b}  F^\col_{b}(\alpha z).
	\end{equation}
\end{proposition}

Note that since $F^\col_b$ is a Coleman function on $U$, by the functorial property of Coleman functions (\cite{F} Proposition 2.5), 
we may regard both $F^\col_b(z+z_\alpha)$ and $F^\col_{b}(\alpha z)$ as Coleman functions on $U_\alpha := E \setminus E[\alpha]$.
Then \eqref{equation: Coleman distribution} above is an equality of Coleman functions on $U_\alpha$.

\begin{proof}[Proof of Proposition \ref{pro: distribution F coleman}]
	We write $\alpha = \pi^n \alpha_0$, where $\alpha_0 \in \cO_\bsK$ is prime to $\pi$.
	For any $z_\alpha \in E[\alpha]$, we write $z_\alpha = z_{\alpha_0} + z_n$, where
	$z_{\alpha_0} \in E[\alpha_0]$ and $z_n \in E[\pi^n]$.
	Suppose $z$ is a point in the residue disc $]z_0[$.  
	Then $z + z_\alpha$ is in the residue disc $]z_0 + z_{\alpha_0}[$, and $\alpha z$ is in $]\alpha z_0[$.  
	Denote by $t$ the parameter on $]z_0[$ obtained as the translation by $\tau_{z_0} : E \rightarrow E$
	of the local parameter $t = -2x/y$ at the origin of $E$.
	This gives an isomorphism $]z_0[ \cong \bbB^-$, hence $t$ parameterizes the points $z$ on the
	residue disc $]z_0[$. 
	We denote  by $t_n$ the
	element on $\bbB^{-}$ corresponding to $z_n$ through this isomorphism.
	Then Proposition \ref{pro: interpolation} shows that we have
	\begin{align*}
		F^\col_{b}(\alpha z)\bigr|_{]z_0[} &= \wh F_{\alpha z_0, b}([\alpha] t),  \\
		F^\col_{b}(z + z_\alpha)\bigr|_{]z_0[}  &= \wh F_{z_0 + z_{\alpha_0}, b}(t \oplus t_n).
	\end{align*}
	as functions on $]z_0[$.
	By Lemma \ref{lem: translation}, we have $\wh F_{z_0+z_{\alpha_0}, b}(t \oplus t_n) = \wh F_{z_0 + z_\alpha, b}(t)$,
	Our result now follows by substituting $z = \lambda(t)$ to the distribution relation of Proposition \ref{pro: distribution F},
	noting that $F_{z_0,b}(\alpha z)|_{z=\lambda(t)} = F_{z_0,b}(z)|_{z=\lambda([\alpha] t)}$.
\end{proof}

\begin{remark}
	\begin{enumerate}
		\item The function $F_1^\col$ is characterized as the unique function of the form $F^\col_1 = \int d F_1$
			satisfying the distribution relation.
		\item The convergence property for $\wh F_{z_0, 1}(t)$ shows that the function $F^\col_1$ 
			converges on any point in $U(\bbC_p)$.
		\item Furthermore, the expansion of $\wh F^\col_b(t)$ for $b>1$ shows that $F^\col_b$ is defined
			on $E(\bbC_p)$.
	\end{enumerate}
\end{remark}

%
\subsection{The $p$-adic Eisenstein-Kronecker function}
%

We will now define the $p$-adic Eisenstein-Kronecker function.  
We first prepare a proposition concerning integration and the distribution relation.

\begin{proposition}\label{pro: distribution}
	Let $m$ and $b$ be integers $\geq 0$.
	Suppose $E_{m,b}^\col$ is a Coleman function on $U$ which satisfies the distribution relation
	\begin{equation}\label{eq: distribution first}
		\sum_{z_\alpha \in E[\alpha]} E^\col_{m,b}(z+z_\alpha) = 
			 \alpha^{1-m} \ol\alpha^{1-b}
		E^\col_{m,b}(\alpha z).
	\end{equation}
	for any non-zero $\alpha \in \cO_\bsK$.  Then there exists a unique integration
	$E^\col_{m+1,b} := -\int E^\col_{m,b} \omega$ of $-E^\col_{m,b}\omega$ satisfying the distribution relation
	\begin{equation}\label{eq: distribution}
		\sum_{z_\alpha \in E[\alpha]} E^\col_{m+1,b}(z+z_\alpha) = \alpha^{-m} \ol\alpha^{1-b}
		E^\col_{m+1,b}(\alpha z) 
	\end{equation}
	for any non-zero $\alpha \in \cO_\bsK$.
\end{proposition}

\begin{proof}
	Let $\wt E_{m+1,b} := -\int E^\col_{m,b} \omega$ be any Coleman integral of 
	$-E^\col_{m,b} \omega$.  For any non-zero $\alpha \in \cO_\bsK$, let
	$$
		c_\alpha: =\sum_{z_\alpha \in E[\alpha]} \wt E_{m+1,b}(z+z_\alpha) - \alpha^{-m} \ol\alpha^{1-b} \wt E_{m+1,b}(\alpha z).
	$$
	Then the relation \eqref{eq: distribution first} shows that $dc_\alpha = 0$, hence the property of Coleman integration
	shows that $c_\alpha$ is a constant in $\bbC_p$.  For any non-zero $\alpha$, $\beta \in \cO_\bsK$, we have
	\begin{multline*}
		\sum_{z_{\alpha\beta} \in E[\alpha\beta]}	\wt E_{m+1,b}(z+z_{\alpha\beta} ) = 
		\sum_{\substack{z_{\alpha} \in E[\alpha], \\\wt z_\beta \in E[\alpha\beta]/E[\alpha]}}
			\wt E_{m+1,b}(z+z_\alpha + \wt z_\beta ) \\
		=	\sum_{\wt z_\beta \in E[\alpha\beta]/E[\alpha]} \left( c_\alpha  +
			\alpha^{-m} \ol\alpha^{1-b}\wt E_{m+1,b}(\alpha z+ \alpha\wt z_\beta )  \right)\\
		=  N(\beta) c_\alpha + \alpha^{-m} \ol\alpha^{1-b} c_\beta
		+  (\alpha\beta)^{-m} (\ol{\alpha\beta})^{1-b} \wt E_{m+1,b}(\alpha \beta z),
	\end{multline*}
	where the last equality follows from the definition of $c_\beta$ and the fact that we have an isomorphism 
	$E[\alpha\beta]/E[\alpha] \cong E[\beta]$ given by $\wt z_\beta \mapsto z_\beta:=\alpha \wt z_\beta$.
	By reversing the roles of $\alpha$ and $\beta$, we see from a similar calculation that the above is also equal to
	$$
		N(\alpha) c_\beta +  \beta^{-m} \ol\beta^{1-b} c_\alpha
		+ (\alpha\beta)^{-m} (\ol{\alpha\beta})^{1-b} \wt E_{m+1,b}(\alpha \beta z).
	$$
	This shows that we have $(N(\beta) - \beta^{-m} \ol\beta^{1-b}) c_\alpha = (N(\alpha) - \alpha^{-m} \ol\alpha^{1-b}) c_\beta$,
	hence the constant
	$$
		c : = c_\alpha/(N(\alpha) - \alpha^{-m} \ol\alpha^{1-b}) 
	$$
	is independent of the choice of $\alpha\in\cO_\bsK$.  If we let 
	$$
		E^\col_{m+1,b}(z) := \wt E_{m+1,b}(z) - c,
	$$
	then this function satisfies \eqref{eq: distribution} for any non-zero $\alpha \in \cO_\bsK$.
\end{proof}

\begin{definition}\label{definition: p-adicEK}
	For integers $m$, $b$ with $b \geq 0$, we define the \textit{$p$-adic Eisenstein-Kronecker series} 
	$E^\col_{m,b}$ recursively going up and down as follows.  We first let
	$E^\col _{0,b}:= (-1)^{b-1} F^\col _{b}$.   Then by Proposition \ref{pro: distribution F coleman}, this function satisfied the
	distribution relation \eqref{eq: distribution first}.
	For $m > 0$,  we let $E^\col_{m,b}$ be the Coleman function recursively defined
	by $E^\col_{m,b} =- \int E^\col_{m-1,b} \omega$, with the constant term normalized as in Proposition \ref{pro: distribution} to
	satisfy the distribution relation
	\begin{equation}\label{eq: distribution final}
		\sum_{z_\alpha \in E[\alpha]} E^\col_{m,b}(z+z_\alpha) =  \alpha^{1-m} \ol\alpha^{1-b} E^\col_{m,b}(\alpha z).
	\end{equation}
	For $m<0$, we define $E^\col_{m,b}$ recursively by the formula $d E^\col_{m+1,b} = -E^\col_{m,b} \omega$.
\end{definition}

Again, \eqref{eq: distribution final} is an equality of Coleman functions on $U_\alpha:= E \setminus E[\alpha]$.
Proposition \ref{pro: distribution} insures that such a choice of constant term
when $m>0$ is possible.  The convergence property of $F_1$ in Proposition \ref{pro: radius convergence} 
insures that $E^\col_{m,b}$ is defined on
any point in $U(\bbC_p)$ if $b=1$ and on $E(\bbC_p)$ if $b > 1$.  When $b=0$, note that $E^\col_{0,0}\equiv 1$.
This shows that we have $E^\col_{a,0} = 0$ for $a<0$.

The reason we view $E^\col_{m,b}$ as a $p$-adic analogue of Eisenstein-Kronecker function is that this function
interpolates values of the classical Eisenstein-Kronecker function at torsion points for $m \leq 0$ as follows.

\begin{proposition}\label{pro: interpolation two}
	Let $a$, $b$ be integers $\geq 0$.  Then for any torsion point $z$ in $E(\ol\bbQ)$ such that $z \not=0$ if $b=1$, we have
	$$
		E^\col_{-a,b}(z) = E_{-a,b}(z)/A^a 
	$$
\end{proposition}

\begin{proof}
	Any torsion point $z$ is of the form $z = z_0 + z_n$, where $z_0$ is a Teichm\"uller representative
	and $z_n$ is a $\pi^n$-torsion point.  Our result follows from the fact that
	$$
		E^\col_{-a,b}(z) = (-1)^{a+b-1} \partial_{t, {\log}}^a \wh F_{z_0,b}(t)|_{t = t_n}
	$$
	and Corollary \ref{cor: generating function}.
\end{proof}

Note that we have fixed a branch of the $p$-adic logarithm.  We next prove that the values of $E^\col_{m,b}(z)$ are 
independent from this choice.

\begin{lemma}\label{lem: IOB}
 	Let $m$ and $b$ be integers such that $b \geq 0$.
	Suppose that $z$ is a point in $E(\bbC_p)\backslash ]0[$.
	Then the value $E_{m,b}^{\mathrm{col}}(z)$ is independent from the choice of
	the branch of the $p$-adic logarithm.
\end{lemma}

\begin{proof}
	The proof is given by induction on $m$.
	The statement for $m=0$ and $b=1$ follows from the construction of $F^\col_1$, which is independent of the
	choice of the branch of the $p$-adic logarithm.  The statement for $m=0$ and general $b$ follows from the 
	definition \eqref{eq: F coleman} of $F_b^\col$, noting that $L_n$ are rational functions.	The statement for $m<0$
	follows inductively from the equality $E^\col_{m,b}  \omega = - d E^\col_{m+1,b}$.
	Suppose now that the statement is true for some $m \geq 0$.
	Let ${\widetilde E}_{m+1,b}=-\int E_{m,b}\omega$ be any Coleman integral.	
	By Besser's formalism [Bes3], Coleman integrals over points of good reduction
	(in other words, points in the smooth subscheme $U \subset X$ using the notations of \S \ref{subsection: Coleman})
	are free from the choice of the branch of the $p$-adic logarithm.
	Hence for $z$ in  $E(\bbC_p)\backslash ]0[$,	the value ${\widetilde E}_{m+1,b}(z)$ is independent of the choice of the branch.
	Furthermore, for a point in $z$ in $E(\bbC_p)\setminus]0[$,
	the values ${\widetilde E}_{m+1,b}(\pi z)$ and ${\widetilde E}_{m+1,b}(z+z_1)$ for $z_1\in E[\pi]$
	are also free from the choice of the branch.
	Hence the global constant $c= c_\pi$ of Proposition \ref{pro: distribution} is independent 
	of the choice of the branch.  This gives the statement for
	$E_{m+1,b}^{\mathrm{col}}(z)$. 
\end{proof}

\begin{remark}\label{rem: IOB}
	The restriction of
	$
		E_{m,b}^{\mathrm{col}}(z)
	$ 
	to the residue disc $]0[$ is of the form
	$$
		E_{m,b}^{\mathrm{col}}|_{]0[}=
		\wh a_0(t)+ \wh a_1(t)\log t+\wh a_2(t)(\log t)^2+\cdots+
		\wh a_n(t)(\log t)^n
	$$
	for some $n$ and $\wh a_i(t) \in A(]0[)$, where the $\wh a_i(t)$ are rigid analytic functions on an open annulus around $0$ free from 
	the choice of the branch.
	Since $E_{m,b}^{\mathrm{col}}(z)$ for $b \neq 1$ is analytic on $E(\bbC_p)$, we see that $n=0$ in this case.
	Therefore if $b \neq 1$, then the values $E_{m,b}^{\mathrm{col}}(z)$ for
	$z\in E(\bbC_p)$ is independent of the choice of the branch.  On the other hand, if $b=1$ and $m>0$, then 
	$E_{m,b}^{\mathrm{col}}(z)$ for $z$ in the residue disc $]0[$ depends on the choice of the branch of the $p$-adic logarithm.
\end{remark}

\section{$p$-adic Kronecker  limit formulas}
%
%
%

Let the notations be as in \S \ref{section: p-adic}.
In this section, we prove the $p$-adic analogues of the first and second Kronecker limit formulas.

%
\subsection{The $p$-adic Eisenstein-Kronecker functions}
%

In \S \ref{section: p-adic}, we defined the $p$-adic analogue of the Eisenstein-Kronecker series as a Coleman function on 
the CM elliptic curve. 
In order to prove the $p$-adic limit formulas, we 
define in this subsection a $p$-adic analogue of the function $\log |\theta(z)|^2-|z|^2/A$, 
which turns out to be a Coleman function.   We then prove the distribution relation, which
we use to characterize this function.

Let $p$ be a prime $\geq 5$.  In what follows, fix an embedding of $\overline{\bbQ}$ into $\bbC_p$ and
we again fix a branch of the $p$-adic logarithm, which is a homomorphism 
$\log_p: \bbC_p^\times \rightarrow \bbC_p$, and we extend this homomorphism 
to $(\bbC_p[[t]])^\times$ by using the decomposition $(\bbC_p[[t]])^\times=\bbC_p^\times \times (1+t\bbC_p[[t]])$ and defining
 $\log_p (1-tf(t))=-\sum {t^nf^n(t)}/n$ for any $f(t) \in \bbC_p[[t]]$. 
Let $E$ be a CM elliptic curve as in the introduction and let $\Gamma$ be the period lattice of $E\otimes\bbC$. 
For $z_0 \in \Gamma \otimes \bbQ$, we let 
 \begin{equation}\label{equation: theta0}
    \theta_{z_0}(z):=\theta(z+z_0) \exp\left(-\frac{z\overline{z_0}}{A}-\frac{z_0\overline{z_0}}{2A}\right).
\end{equation}
Then by \cite{BK1} Proposition 2.1, the Taylor  series of $\theta_{z_0}(z)$ at $z=0$ has algebraic coefficients.
Note that from the definition, we have 
$$	\frac{d}{dz} \log \theta_{z_0}(z)=
	\frac{\theta'_{z_0}(z)}{\theta_{z_0}(z)} = \frac{\theta'(z+z_0)}{\theta(z+z_0)} - \frac{\ol z_0}{A}
	= F_1(z+z_0)  - \frac{\ol z_0}{A} = F_{z_0,1}(z).
$$
If we consider the formal composition 
$$
    \widehat{\theta}_{z_0}(t):=\theta_{z_0}(z)|_{z=\lambda(t)}
$$
of this series with $\lambda(t)$, where $\lambda(t)$ is the formal logarithm of the formal group of $E$, then we may regard this 
power series as an element in  $\bbC_p[[t]]$. 
If we take the derivative of $\log_p \widehat{\theta}_{z_0}(t)$ with respect to $t$, then we have by the definition of
$\log_p$ and $\wh\theta_{z_0}(t)$ the equality
\begin{equation}\label{equation: derivative}
	\bigl( \log_p \widehat{\theta}_{z_0}(t) \bigr)' = \widehat{\theta}'_{z_0}(t) / \widehat{\theta}_{z_0}(t) = \wh F_{z_0, 1}(t).
\end{equation}
If $z_0$ corresponds to a torsion point of $E$ of order prime to $\frp$, then 
we see from Proposition \ref{pro: radius convergence} that the the power series
$\bigl( \log_p \widehat{\theta}_{z_0}(t) \bigr)' $ hence also 
$\log_p \widehat{\theta}_{z_0}(t)$ converges on the open unit disc $|t|<1$ in $\bbC_p$. 


\begin{definition}\label{def: log-p-theta}
	We let $\log _p \theta$ be the function in $A_{\loc}(U)$ defined by 
	$$
    		\log _p \theta |_{]z_0[}:=\log_p \widehat{\theta}_{z_0}(t) \in A_{\log}(]z_0[)
	$$
	on each residue disc $]z_0[$, where $z_0$ corresponds to a torsion point of $E$ of order prime to $\frp$.
\end{definition}

We will see in \S \ref{subsection: 5.2} that $\log_p {\theta}$ is in fact a Coleman function on $E$.
We first investigate the basic properties of $\log_p {\theta}$. 

\begin{proposition} 
	For $z_0 \in \Gamma \otimes \bbQ$ and $z_\alpha$ such that 
	$\alpha z_\alpha \in \Gamma$ for a $\pi$-power morphism $\alpha \in \mathrm{End}_{\overline{\bbQ}}(E)$, we have 
	\begin{equation}\label{equation: proposition}
  		\log_p \widehat{\theta}_{z_0}(t \oplus t_\alpha)=\log_p\widehat{\theta}_{z_0+z_\alpha}(t),
 	\end{equation}
	where $t_\alpha \in E(\bbC_p)$ corresponds to a torsion point
	in $\bbC/\Gamma$ represented by $z_\alpha \in \bbC$.
\end{proposition}

\begin{proof}
          Let $\alpha$ and $\beta$ be elements of $\cO_K$ such that 
          $2\alpha | \beta$ and $\beta z_0 \in \Gamma$. 
           Then  $f_{\beta}(z):=\theta(z)^{N\beta}/\theta( \beta z)$  is a rational function on $E$ over $\overline{\bbQ}$.  
           We have 
 	\begin{equation}\label{equation: logtheta}
    		\theta_{z_0+z_\alpha}(z)^{N\beta}= \pm \theta( \beta z) 
 		\tau_{z_0+z_\alpha}^*f_{\beta}(z). 
 	\end{equation}
	Similarly,  we have 
  	$$
        		\theta_{z_0}(z)^{N\beta}
                  = \pm \theta( \beta z) \tau_{z_0}^*f_{\beta}(z). 
  	$$
	Since  $f_\beta$ is a rational function, we also have 
  	$$
      		\tau_{z_0+z_\alpha}^*f_{\beta}(t)=  \tau_{z_0}^*f_{\beta}(t\oplus t_\alpha).
  	$$
	Hence  we have 
   	\begin{equation}\label{equation: logtheta2}
    		\widehat{\theta}_{z_0}(t \oplus t_\alpha)^{N\beta}
                  = \pm\widehat{\theta}( [\beta] t) \tau_{z_0+z_\alpha}^*f_{\beta}(t). 
  	\end{equation}
	Our assertion now follows from \eqref{equation: logtheta} and \eqref{equation: logtheta2}.
\end{proof}

\begin{corollary}\label{corollary: interpolation}
	Let $t_\alpha$ be a $\pi$-power torsion point, and we assume that $z_0 \not=0$ or $t_\alpha\not=0$. 
	Then we have 
	$$
  		\log_p \widehat{\theta}_{z_0}(t_\alpha)=
		\log_p \left(\theta(z_0+z_\alpha) \exp\left[-\frac{(z_0+z_\alpha)\overline{(z_0+z_\alpha)}}{2A}\right] \right).
	$$
\end{corollary}

\begin{proof}
	This follows by substituting $t=0$ on both sides of \eqref{equation: proposition} and using
	the definition of $\theta_{z_0 + z_\alpha}(z)$ given in \eqref{equation: theta0}.
\end{proof}

Roughly speaking,  $\log_p \theta(z)$ is a $p$-adic function which interpolates the special values $\log \theta(z)-z \overline{z}/2A$ at torsion points. 
We may thus regard  $\log_p \theta(z)$ as a $p$-adic analogue of the function $\log |\theta(z)|^2-|z|^2/A$.

%
\subsection{The $p$-adic second limit formula}\label{subsection: 5.2}
%

We are now ready to prove Theorem \ref{thm: pKLF}, which is a $p$-adic analogue of the Kronecker
second limit formula.  We keep the notation of the introduction.  In particular, let
$$
	K^\col_{n-m}(0,z,n) := E^\col_{m,n}(z).
$$
In addition, we let $\frf \subset \cO_K$ be the conductor of the Hecke character $\psi$ of $K$ associated to $E$,
and $\pi := \psi(\frp)$.
Since $E$ has good reduction at the prime ideals above $p$, we have $(p, \frf)=1$.
\begin{proposition}\label{proposition: theta distribution}
         Let $z_0 \in \bbC$ be a lifting of a $\mathfrak{f}$-torsion point of $\bbC/\Gamma$. Then for $\alpha \in \cO_K$ 
	we have 
	$$
      	\theta_{\alpha z_0}( \alpha z)^{24 N(\alpha\mathfrak{f})}= \Delta^{2N(\alpha\mathfrak{f})(N\alpha-1)} 
                        \prod_{z_\alpha \in E[\alpha]} \theta_{z_0+z_\alpha}(z)^{24 N(\alpha\mathfrak{f})},
	$$
where $z_\alpha$ is a lift of a $\alpha$-torsion point of $E$ and 
the right hand side is independent of the choice of the lifts $z_0$ and $z_\alpha$ on $\bbC$. 
\end{proposition}

\begin{proof}
         Since for $\gamma \in \Gamma$  we have 
         $\theta_{z_0+\gamma}(z)=\pm \chi_\gamma(z_0/2) \theta_{z_0}(z)$, 
         the function $ \theta_{z_0}(z)^{2N\mathfrak{f}}$ is independent of the lift $z_0$ if 
         $(N\mathfrak{f})z_0 \in \Gamma$. The independence of the lifts of $z_0$ and $z_\alpha$ follows from this fact. 
	The logarithmic  derivatives of both sides coincide by Proposition \ref{pro: distribution F}. 
	Hence for each $\alpha$,  there exists a constant $c_\alpha(z_0)$ such that 
	 $$
     		\theta_{\alpha z_0}( \alpha z)^{2N(\alpha\mathfrak{f})}= c_\alpha(z_0)
                        \prod_{z_\alpha \in E[\alpha]} \theta_{z_0+z_\alpha}(z)^{2N(\alpha\mathfrak{f})}. 
	 $$
	 If we compare this equality with the case $z_0 = 0$, then we have from the definition of $\theta_{z_0}(z)$
	 given in \eqref{equation: theta0}
	 the equality
	 $$
	 	c_\alpha(0) = c_\alpha(z_0) \prod_{z_\alpha} \exp\left( \frac{z_0 \ol z_\alpha - \ol z_0 z_\alpha}{2A}\right)^{2N(\alpha \frf)}
		= c_\alpha(z_0).
	 $$
	Hence we see that $c_{\alpha} := c_\alpha(z_0)$ is independent of  the choice of $z_0$. 
	We calculate $c_\alpha$ for the case $z_0=0$ and $N\mathfrak{f}=1$. 
	Then we have 
	\begin{align*}
		 \prod_{z_{\alpha \beta} \in E[\alpha \beta]} \theta_{z_{\alpha \beta}}(z)^{2N(\alpha\beta)}
		= \prod_{z_{\alpha \beta} \in E[\alpha \beta]/E[\alpha]} \prod_{z_{\alpha } \in E[\alpha ]} 
		  \theta_{z_{\alpha \beta}+z_{\alpha}}(z)^{2N(\alpha\beta)}\\
		= \prod_{z_{\alpha \beta} \in E[\alpha \beta]/E[\alpha]}
		c_\alpha^{-N\beta}  \theta_{ \alpha z_{\alpha \beta}}(\alpha z)^{2N(\alpha\beta)}\\
		= c_\alpha^{-N\beta^2} c_\beta^{-N\alpha} \theta(\beta \alpha z)^{2N(\alpha\beta)}.
	\end{align*} 
	Hence we have 
 	$
   		c_\alpha^{N\beta^2} c_\beta^{N\alpha}=c_{\alpha \beta}=c_\beta^{N\alpha^2} c_\alpha^{N\beta}
 	$
	or equivalently,  
	$$
   		c_\alpha^{N\beta(N\beta-1)}=c_\beta^{N\alpha(N\alpha-1)}. 
 	$$
	In particular, $c_\alpha^{12}=c_2^{N\alpha(N\alpha-1)}$. 
	On the other hand, we consider the constant term of 
	$$
   		\frac{\theta(2z)^8}{\theta(z)^8}=c_2   \prod_{z_2 \in E[2]-\{0\}} \theta_{z_2}(z)^{8}. 
	$$
	As in the proof of Proposition \ref{pro: C}, 
	we have 
	$$
  		\prod_{z_2 \in E[2]-\{0\}} \theta_{z_2}(0)^{8}=\Delta'^{-2}.
	$$
	Hence $c_2=2^8 \Delta'^2=\Delta^2$. Our assertion now follows from these facts. 
\end{proof}

\begin{corollary}\label{corollary: 5.6}
	The function $\Xi(z):= - \log_p \theta (z) -\frac{1}{12}\log_p \Delta$ satisfies the distribution relation  
	$$
		\Xi(\alpha z)=\sum_{z_\alpha \in E[\alpha]} \Xi(z+z_\alpha). 
	$$
\end{corollary}

\begin{proof}
	By Proposition \ref{proposition: theta distribution}, we have 
	$$
		\log_p \widehat{\theta}_{\alpha z_0}([\alpha]t)=\frac{N\alpha-1}{12} \log_p \Delta 
		+\sum_{t_\alpha \in E[\alpha]} \log_p \widehat{\theta}_{z_0}(t \oplus t_\alpha)
	$$
	on each residue disc $]z_0[$. 
	Our assertion follows from this formula, since 
	\begin{multline*}
		\Xi(\alpha z)\bigr|_{]z_0[} :=- \log_p \wh\theta_{\alpha z_0} ([\alpha]t) -\frac{1}{12}\log_p \Delta\\
		=-\sum_{t_\alpha \in E[\alpha]} \log_p \widehat{\theta}_{z_0}(t \oplus t_\alpha)- \frac{N\alpha}{12} \log_p \Delta 
		= \sum_{z_\alpha \in E[\alpha]} \Xi(z+z_\alpha)\bigr|_{]z_0[}
	\end{multline*}
	on each residue disc.
\end{proof}

We now prove the $p$-adic second limit formula.

\begin{proof}[Proof of Theorem \ref{thm: pKLF}]
	By the definition of $\Xi(z)$ and \eqref{equation: derivative}, the derivative of  $\Xi(z)$ is equal to 
	$-\wh F_{z_0, 1}(t)$ on the residue disc $]z_0[$.  By Definition \ref{definition: p-adicEK}, the derivative of
	$$
		K_0^\col(0,z,1):=E^{\col}_{1,1}(z)
	$$ 
	is equal to $E^\col_{0,1}(z) := -F_1^\col(z)$, hence coincides with $-\wh F_{z_0, 1}(t)$ on $]z_0[$.
	This implies that $c(z):=\Xi(z)-E^{\col}_{1,1}(z)$ is a constant on  the residue disc $]z_0[$. 
	By Corollary \ref{corollary: 5.6} and the definition of $E^{\col}_{1,1}(z)$, the locally constant function $c(z)$ satisfies the distribution relation. 
	For any torsion point $z_0$ of order $\mathfrak{f}$, we take $N$ such that $\pi^N \equiv 1 \mod \mathfrak{f}$. 
	Then $[\pi^N]^*(]z_0[)=]z_0[$ and 
	$$
		[\pi^N]^*c(z) |_{]z_0[}=\sum_{w \in E[\pi^N]} c(z+w) |_{]z_0[}. 
	$$
	Since $c(z) |_{]z_0[}$ is constant, the above relation shows $c(z) |_{]z_0[}=0$.  
\end{proof}

The above result shows in particular that  
$$
	\Xi(z)=- \log_p \theta (z) -\frac{1}{12}\log_p \Delta
$$ 
is in fact a Coleman function.

%
\subsection{$p$-adic Eisenstein-Kronecker-Lerch series}
%

We now give the definition of the $p$-adic Eisenstein-Kronecker-Lerch series and then
prove Proposition \ref{pro: pKLF}, which is a $p$-adic analogue of the first Kronecker limit formula.
We will prove the proposition by considering $p$-adic counterparts of our proof in \S \ref{section: one} of the classical case.

Let $p \geq 5$ be a prime of good \textit{ordinary} reduction for $E$, and we fix a prime $\frp$ of $\cO_\bsK$ over $p$.  
We defined in \cite{BK1} \S 3.1 a $p$-adic measure $\mu:=\mu_{0, 0}$ on $\bbZ_p \times \bbZ_p$ interpolating
the Eisenstein-Kronecker numbers, or more precisely,
the special values of Eisenstein-Kronecker-Lerch series $K^*_{a+b}(0,0,b; \Gamma)/A(\Gamma)^a$ for $a$, $b \geq 0$,
where $\Gamma$ is the period lattice of $E$.
We define the $p$-adic Eisenstein-Kronecker-Lerch function as in the introduction as follows.

\begin{definition}  For any integer $a \in \bbZ$, we define the \textit{$p$-adic Eisenstein-Kronecker-Lerch function} by
	$$
		K^{(p)}_a(0,0,s) := \int_{\bbZ_p^\times \times \bbZ_p^\times} \pair{x}^{s-1}\pair{y}^{a-s} \omega(y)^{a-1} d \mu(x,y).
	$$
\end{definition}


The $p$-adic Eisenstein-Kronecker-Lerch function is analytic in $s \in \bbZ_p$. 
The reason we view this function as a $p$-adic analogue of Eisenstein-Kronecker-Lerch series is the
following interpolation property.

\begin{proposition}
	For any integer $a$, $b$ such that $a \geq b > 0$ and $b \equiv 1 \pmod{p-1}$, we have
	\begin{equation}\label{eq: interpolation}
		\frac{K^{(p)}_a(0,0,b)}{\Omega_p^{a-1}} 
		=(-1)^{a-1}(b-1)! \left(1-\frac{\pi^{a}}{p^{a-b+1}}\right)
		\left(1-\frac{\pi^{a}}{p^b} \right)\frac{K^*_{a}(0,0,b)}{A(\Gamma)^{a-b}},
	\end{equation}
	where $\Omega_p$ is a $p$-adic period of the formal group of $E$.
\end{proposition}

\begin{proof}
	This follows from the interpolation property of the measure $\mu:=\mu_{0,0}$ given in \cite{BK1} Proposition 3.5.
\end{proof}

We now give the proof of Proposition \ref{pro: pKLF}.  

\begin{proof}[Proof of Proposition \ref{pro: pKLF}]
	We consider the function 
	\begin{equation}\label{eq: f}
		f(t):=\Omega_p \int_{\bbZ_p^\times \times \bbZ_p^\times}y^{-1} \exp(y\Omega_p^{-1} \lambda(t)) d\mu(x,y)
	\end{equation}
	on the $p$-adic residue disc $]0[$ around $0$.
	If we take the derivative of $f(t)$, the interpolation property of $\mu$ given in \cite{BK1} Proposition 3.5 gives the equality
	\begin{multline*}
 		\lambda'(t)^{-1}\frac{d}{dt}f(t)=\int_{\bbZ_p^\times \times \bbZ_p^\times} \exp(y\Omega_p^{-1} \lambda(t)) d\mu(x,y)\\
		= \wh F_1(t; \Gamma)-{\overline{\pi}}^{-1}\wh F_1([\pi]t; \Gamma)-\wh F_1(t; \ol\frp\Gamma)
		 +{\overline{\pi}}^{-1} \wh F_1([\pi]t; \ol\frp\Gamma).
	\end{multline*}
	Let $E_{1,1}^{(p)}(z; \Gamma) := E_{1,1}^\col(z; \Gamma) - p^{-1} E_{1,1}^\col(\pi z; \Gamma)$.  Then the differential of 
	$
		E_{1,1}^{(p)}(z; \Gamma)-E_{1,1}^{(p)}(z; \ol\frp\Gamma)
	$
	is given by
	$$
		F^\col_1(z; \Gamma)\omega - \frac{\pi}{p} F^\col_1(\pi z;\Gamma) \omega
		+ F^\col_1(z; \ol\frp\Gamma)\omega - \frac{\pi}{p} F^\col_1(\pi z;\ol\frp\Gamma) \omega,
	$$
	which is equal to $d f(t)$ on $]0[$, hence the function $f(t) - E_{1,1}^{(p)}(z; \Gamma)+E_{1,1}^{(p)}(z; \ol\frp\Gamma)$
	is a constant on the residue disc $]0[$. 
	By the definition of $E_{1,1}^{(p)}$ and substituting $z = \lambda(t)$ into \eqref{eq: distribution}, we have
	the distribution relation  
	$$\sum_{t_{\pi} \in E[\pi]} \wh E_{1,1}^{(p)}(t\oplus t_{\pi})=0.$$ 
	Furthermore, since 
	the power series $\exp(\Omega_p^{-1} \lambda(t))$ gives a homomorphism of
	formal groups $\wh E$ and $\wh \bbG_m$ isomorphically mapping $E[\pi]$ to the group of $p$-th root of unity
	\cite{BKT1} \S 2.2, we have for each $y \in \bbZ_p^\times$ the equality
	$$
		\sum_{t_\pi \in E[\pi]}\exp( y\Omega_p^{-1} \lambda(t))\big|_{t=t \oplus t_\pi} =
			\exp( y\Omega_p^{-1} \lambda(t))  \sum_{t_\pi \in E[\pi]}\exp( y\Omega_p^{-1} \lambda(t_\pi))  = 0.
	$$
	This gives by definition \eqref{eq: f} of $f(t)$ the distribution relation $\sum_{t_{\pi} \in E[\pi]} f(t\oplus t_{\pi})=0$.
	Since both $ E_{1,1}^{(p)}$ and $f(t)$ satisfy the same distribution relation, the same argument as that in the proof of Theorem \ref{thm: pKLF}
	shows that we must have 
	$f(t)=E_{1,1}^{(p)}(z; \Gamma)-E_{1,1}^{(p)}(z; \ol\frp\Gamma)$ 
	on ${]0[}$. On the other hand, the $p$-adic  second limit formula shows that 
	$$
		E_{1,1}^{(p)}(z; \Gamma)=\log_p \theta(z; \Gamma)
		-\frac{1}{p}\log_p \theta(\pi z; \Gamma)+\frac{1}{12}\left(1-\frac{1}{p}\right) \log_p \Delta(\Gamma).
	$$
	Noting that 
	$$
		\log_p \theta(z; \Gamma) -\log_p \theta(z; \ol\frp\Gamma) = \log_p \left(\theta(z; \Gamma)/\theta( z; \ol\frp\Gamma) \right) 
	$$
	for $z=0$ is equal to $\log_p(1) = 0$ and $\Delta(\ol\frp \Gamma) = \ol\pi^{-12} \Delta(\Gamma)$, we have
	$$
		f(0) = \bigl.E_{1,1}^{(p)}(z; \Gamma)-E_{1,1}^{(p)}( z; \ol\frp\Gamma)\;\bigr|_{z=0}=
		\left(1-\frac{1}{p}\right) \log_p \overline{\pi}.
	$$
	Our assertion now follows from the fact that
	$
		f(0)=\Omega_p K^{(p)}_0(0,0,1).
	$
\end{proof}

\begin{remark}\label{rem: not pKLF}
	In the interpolation formula of \eqref{eq: interpolation}, if we let $a=0$ and $b=1$,  then
	the interpolation factor of the the right hand side vanishes. Hence the value 
	$$
		\Omega_p K^{(p)}_0(0,0,1) = \int_{\bbZ_p^\times \times \bbZ_p^\times}  y^{-1} d\mu(x,y)
	$$
	is in some sense not the constant term $0$ but 
	the residue at $s=1$ of the $p$-adic analogue of $\sum_{\gamma \in \Gamma}^*1/|\gamma|^{2s}$. 
	Because of this fact, the formula of Proposition \ref{pro: pKLF} is not a perfect $p$-adic analogue of the classical first Kronecker limit formula. 
\end{remark}





\begin{thebibliography}{999}
	\bibitem[Ba1]{Ba1}{K. Bannai},  Rigid syntomic cohomology and $p$-adic polylogarithms, J. Reine Angew. Math., 
	\textbf{529} (2000), 205-237.
	\bibitem[Ba2]{Ba3}{K. Bannai}, On the $p$-adic realization of elliptic polylogarithms for CM-elliptic curves, 
	Duke Math. J., \textbf{113} (2002), no. 2, 193-236.
	\bibitem[Ba3]{Ba4}{K. Bannai}, Specialization of the $p$-adic polylogarithm to $p$-th power roots of unity,
	Doc. Math. (2003), Extra Vol., 73--97.
	\bibitem[BK1]{BK2}{K. Bannai and S. Kobayashi}, Algebraic theta functions and
	Eisenstein-Kronecker numbers, RIMS K\^oky\^uroku Bessatsu \textbf{B4}: \textit{Proceedings of the Symposium on 
	Algebraic Number theory and Related Topics},  eds. K. Hashimoto, Y. Nakajima and H. Tsunogai, December (2007), 63--78.
	\bibitem[BK2]{BK1}{K. Bannai and S. Kobayashi}, Algebraic theta functions and 
	$p$-adic interpolation of Eisenstein-Kronecker numbers,  Duke Math. J. \textbf{153} no. 2 (2010), 229-295.
	\bibitem[BKT1]{BKT2}{K. Bannai, S. Kobayashi, and T. Tsuji}, Realizations of the elliptic polylogarithm for CM elliptic curves,
	RIMS K\^oky\^uroku Bessatsu \textbf{B12}: \textit{Algebraic Number Theory and Related Topics 2007}, 
	eds. M. Asada, H. Nakamura and H. Takahashi, August (2009), 33--50.
	\bibitem[BKT2]{BKT1}K. Bannai, S. Kobayashi and T. Tsuji, 
	On the de Rham and $p$-adic realizations of the elliptic polylogarithm 
	for CM elliptic curves, \textit{Annales scientifiques de l'ENS} \textbf{43}, fascicule 2 (2010), 185-234.
	\bibitem[BL]{BL}{A. Beilinson and A. Levin}, The Elliptic Polylogarithm, 
	\textit{Motives}, (Seattle, WA, 1991), 123-192. 
	\bibitem[Ber]{Ber}{P. Berthelot}, Finitude et puret\'e cohomologique en cohomologie rigide,
	avec un appendice de A. J. de Jong, Inventiones Math. \textbf{128},  329-377 (1997).	
	\bibitem[Bes1]{Bes1}{A. Besser}, Syntomic regulators and p-adic integration I: Rigid syntomic regulators, 
		Israel Journal of Math. \textbf{120} (2000), 291-334.
	\bibitem[Bes2]{Bes2}{A. Besser}, Syntomic regulators and p-adic integration II: $K_2$ of curves, 
		Israel Journal of Math. \textbf{120} (2000), 335-360.
	\bibitem[Bes3]{Bes3}{A. Besser}, Coleman integration using the Tannakian formalism,
	Mathematische Annalen \textbf{322} (2002) 1, 19-48.
	\bibitem[BdJ]{BdJ}{A. Besser and R. de Jeu}, The syntomic regulator for the $K$-theory of fields.  Ann. Sci. \'Ecole 
	Norm. Sup. (4)  \textbf{36}  (2003),  no. 6, 867--924.
	\bibitem[C1]{Col1} {R. Coleman}, Dilogarithms, regulators and $p$-adic $L$-functions,
	Invent. Math. \textbf{69} (1982), no. 2, 171--208. 
	\bibitem[C2]{Col2} {R. Coleman}, Torsion points on curves and $p$-adic abelian integrals,
	Ann. of Math. (2) \textbf{121} (1985), no. 1, 111--168. 
	\bibitem[D]{D}{P. Deligne}, Le groupe fondamental de la droite projective moins trois points, In:
	\textit{Galois groups over $\bbQ$} (Berkeley, CA, 1987), 79--297, 
	\bibitem[dS]{dS}{E. de Shalit}, \textit{Iwasawa theory of elliptic curves with complex multiplication},
	Academic Press (1987).  
	\bibitem[F]{F} H. Furusho, $p$-adic multiple zeta values I., Inv. Math. \textbf{55} (2004), 253-286.
	\bibitem[K]{Ka}N. Katz, $p$-adic interpolation of real analytic Eisenstein
	series, Ann.\ of Math.\ \textbf{104} (1976), 459-571.
	\bibitem[L]{L}{A. Levin}, Elliptic polylogarithm: An analytic theory, Comp. Math. \textbf{106} (1997), 267-282. 
	\bibitem[R]{Rob} G. Robert, Unit\'es elliptiques et formules pour le nombre de classes des extensions ab\'eliennes
	d'un corps quadratique imaginaire, Bull. Soc. Math. France M\'emoire \textbf{36} (1973).
	\bibitem[S]{Som}{M. Somekawa}, Log-syntomic regulators and $p$-adic polylogarithms,
	 $K$-Theory  \textbf{17}  (1999),  no. 3, 265--294.
        \bibitem[We]{We}{A. Weil}, {\it Elliptic Functions according to Eisenstein and Kronecker},
         Springer-Verlag, 1976. 
	\bibitem[Wi]{W2}{J. Wildeshaus}, On an elliptic analogue of Zagier's conjecture,  Duke Math. J.  \textbf{87} (1997),  
	no. 2, 355--407. 
	\bibitem[Z]{Zag}{D. Zagier}, Periods of modular forms and Jacobi theta functions, Inv. math. 
	\textbf{104} (1991), 449-465.
\end{thebibliography}
\end{document}